\documentclass[11pt,a4paper]{amsart}
\usepackage[english]{babel}
\usepackage[latin1]{inputenc}
\usepackage{color}
\usepackage{graphics}

\usepackage{amsmath,amssymb,amstext,amsthm,amscd}

\def\ligne{\begin{center}
          \unitlength=1cm
          \begin{picture}(5,0.05)
          \line(1,0){5}
          \end{picture}
          \end{center}}

\newcounter{cimage}
\def\C{\mathbb C }
\def\R{\mathbb R }

\def\G{\mathbb G }
\def\F{\mathbb F }
\def\P{\mathbb P }

\def\k{\mathfrak k }
\def\g{\mathfrak g }
\def\z{\mathfrak z }
\def\S{\mathcal S}
\def\V{\mathcal V}

\def\deg{\mathrm{deg}}

\def\tr{\mathrm{tr}}

\def\Herm{\mathrm{Herm}}

\def\leq{\leqslant}

\def\qmod#1#2{{\hbox{}^{\displaystyle{#1}}}\!\big/\!\hbox{}_{
\displaystyle{#2}}}

\newfam\msbfam

\font\tenmsb=msbm10
\font\sevenmsb=msbm10 at 7pt
\font\fivemsb=msbm10 at 5pt

\textfont\msbfam=\tenmsb
\scriptfont\msbfam=\sevenmsb
\scriptscriptfont\msbfam=\fivemsb

\def\map{\longrightarrow}
\def\textmap#1{\mathop{\vbox{\ialign{
                                  ##\crcr
      ${\scriptstyle\hfil\;\;#1\;\;\hfil}$\crcr
      \noalign{\kern-1pt\nointerlineskip}
      \rightarrowfill\crcr}}\;}}

\def\textlmap#1{\mathop{\vbox{\ialign{
                                  ##\crcr
      ${\scriptstyle\hfil\;\;#1\;\;\hfil}$\crcr
      \noalign{\kern-1pt\nointerlineskip}
      \leftarrowfill\crcr}}\;}}

\newfam\meuffam

\font\tenmeuf=eufm10
\font\sevenmeuf=eufm7

\textfont\meuffam=\tenmeuf
\scriptfont\meuffam=\tenmeuf
\scriptscriptfont\meuffam=\sevenmeuf

\def\germ{\fam\meuffam\tenmeuf}

\def\g{{\germ g}}

\def\kg{{\germ k}}

\def\EE{{\cal E}}
\def\FF{{\cal F}}

\def\tr{{\rm Tr}}
\def\End{{\rm End}}
\def\Aut{{\rm Aut}}

\def\U{{\rm U}}

\def\ub{\underbar}

\def\deg{{\rm deg}}
\def\Hom{{\rm Hom}}
\def\Aut{{\rm Aut}}

\def\Herm{{\rm Herm}}
\def\Vol{{\rm Vol}}

\def\id{{\rm id}}

\def\im{{\rm im}}
\def\rk{{\rm rk}}

\def\ad{{\rm ad}}

\def\Ad{{\rm Ad}}

\def\ad{{\rm ad}}

\def\cal{\mathcal}

\def\ra{\rightarrow}

\title{Harder-Narasimhan filtrations and optimal destabilizing
vectors in complex geometry}

\author{L. Bruasse$^{\dagger,1}$ and A. Teleman$^\ddagger$}

\begin{document}
\newtheorem{toto}{toto}[section]
\def\thms{toto}


\theoremstyle{definition}
\newtheorem{definition}[\thms]{Definition}
\newtheorem{example}[\thms]{Example}
\newtheorem*{notations}{Notations}

\theoremstyle{plain}
\newtheorem{theorem}[\thms]{Theorem}
\newtheorem*{theorem2}{Theorem}
\newtheorem{proposition}[\thms]{Proposition}
\newtheorem{property}[\thms]{Property}
\newtheorem{lemma}[\thms]{Lemma}
\newtheorem{corollary}[\thms]{Corollary}

\theoremstyle{remark}
\newtheorem{remark}[\thms]{Remark}
\newtheorem*{notation}{Notation}
\theoremstyle{remark}
\newtheorem{exemple}[\thms]{Example}
\maketitle

\begin{center}
\small \it $^\dagger$IML, CNRS UPR 9016 \\ 163 avenue de Luminy,
case 930, 13288 Marseille cedex 09 - France \\
$^\ddagger$CMI, LATP UMR 6632 \\ 39, rue F. Joliot Curie, 13453
Marseille cedex 13 - France\textcolor{white}{\footnotemark}
\end{center}

\footnotetext{Email addresses: bruasse@iml.univ-mrs.fr (L.
Bruasse), teleman@cmi.univ-mrs.fr (A. Teleman)}

\ligne \vspace*{.2cm}
\begin{center}
\parbox{12cm}{%
\small \textsc{{Abstract.}} We give a
generalisation of the theory of optimal destabilizing 1-parameter
subgroups to non-algebraic complex geometry. Consider  a holomorphic action
$G\times F\to F$ of a complex reductive Lie group $G$ on a finite dimensional
(possibly non-compact) K\"ahler manifold $F$. Using a Hilbert type criterion
for the (semi)stability of symplectic actions, we associate to any non
semistable point $f\in F$ a unique optimal destabilizing vector in $\g$ and
then a naturally defined point
$f_0$ which is semistable for the action of a certain reductive subgroup of $G$
on a submanifold of $F$. We get a natural stratification of $F$
which is the analogue of the Shatz stratification for holomorphic
vector bundles. In the last chapter we  show that our
results can be generalized to the gauge theoretical framework: first 
we show that
the system of semistable quotients associated with the classical 
Harder-Narasimhan filtration of a non-semistable bundle $\EE$
can be recovered as the limit object in the direction given by the 
optimal destabilizing
vector of $\EE$. Second, we extend this principle to holomorphic 
pairs: we give the analogue of the Harder-Narasimhan theorem  for 
this moduli problem and we discuss the  relation between the 
Harder-Narasimhan filtration of a non-semistable holomorphic pair and 
its optimal destabilizing vector.

\

}
\end{center}\vspace*{.2cm}

\ligne


\setcounter{section}{0}
\section{Introduction}

A classical result of Harder and Narasimhan states that any non-semi\-stable
bundle on a curve admits a canonical filtration of subsheaves with torsion free
semistable quotients.

This result was  generalized for reflexive sheaves on projective
varieties \cite{sha}, \cite{mar}, and finally to reflexive sheaves
on arbitrary compact Hermitian  manifolds \cite{bru1}, \cite{bru2}.

The initial motivation for this paper was to find the analogous statement for
other type of complex geometric objects, for instance holomorphic bundles
coupled with sections or with endomorphisms (Higgs fields).

The system  of semistable quotients associated with the
Harder-Narasim\-han filtration of a non-semistable bundle can be
interpreted as a  semistable object with respect to the moduli
problem for $G$-bundles, where $G$ is a product of reductive group
of the form $\prod_i GL(r_i)$.

Therefore, the Harder-Nara\-sim\-han result can be
understood as an assignment which associates to a non-semistable object a
semistable object but for a different moduli problem.

We believe that it is a natural and important problem to find a
general principle which generalizes  this result for
\ub{arbitrary} moduli problems. More precisely, we seek a general
rule which  associates -- in a canonical way -- to a
non-semistable object ${\cal O}$ with respect to \ub{any} complex
geometric moduli problem ${\cal A}$ a new moduli problem ${\cal
B}({\cal A},{\cal O})$ and a semistable object ${\cal O}'({\cal
A},{\cal O})$  for  ${\cal B}({\cal A},{\cal O})$.

Our first attempt was to understand this principle in the finite
dimensional framework, i.e. for   moduli problems associated with
actions of reductive groups on finite dimensional but in general {\it 
non-compact} varieties.

   After consulting
the available literature dedicated to the algebraic case, we
realized that the main tool for understanding the analogue of the
Harder-Narasimhan assignment in the finite dimensional framework
is {\it the theory of optimal one-parameter subgroups}, for which we 
refer to Kirwan  \cite{ki}, Slodowy
\cite{slo} and Ramanan\&Ramanathan \cite{ram}.

This theory can be sketched as follows: if
$[x]\in\P^n(V)$ is non-semistable point with respect to a linear
representation $\rho:G\ra GL(V)$ of a reductive group $G$, then
there exists a one parameter subgroup $\tau:\C^*\ra G$ of "norm" 1
which "destabilizes" $[x]$ in the strongest possible way, i. e.
$$\lambda(x,\tau)\leq \lambda(x,\theta),
$$
for any one-parameter subgroup   $\theta:\C^*\ra G$ of norm 1. Here
we denoted by $\lambda$ the maximal weight function which occurs in
the Hilbert criterion for stability.   A one parameter subgroup
    (an OPS) with this property is called an  {\it optimal
destabilizing} OPS for
$[x]$ and it is essentially unique, in the sense that any other {
optimal} destabilizing OPS
$\tau'$ for $[x]$ has the same associated parabolic subgroup as
$\tau$, and is conjugated with $\tau$ in this parabolic subgroup.

This result has certainly become part of classical GIT. What is
(at least for the authors) less standard material is the  following
crucial remark, which was probably first pointed out by 
Ramanan\&Ramanathan \cite{ram} (see also Kirwan \cite{ki}):\\

{\it If $\tau:\C^*\ra G$ is an  {\it optimal} destabilizing
OPS for $[x]$, then $\tau(t)$ converges to  a point $[x_0]$ which is
semistable with respect to an induced action of the reductive centralizator
$Z(\tau)$ of $\tau$ on a  $Z(\tau)$-stable subvariety of $\P(V)$. }\\

We claim that the assignment $[x]\mapsto [x_0]$ is the GIT model
which should be followed in order to get the correct
generalization of the Harder-Narasimhan theorem to the non-compact
non-algebraic and gauge theoretical frameworks.

Therefore our final goal is to give a gauge theoretical version of
this remark, which applies to \ub{all} moduli problems obtained by
coupling holomorphic bundles (with arbitrary reductive structure
groups) with sections in associated bundles (see \cite{mu},
\cite{lute}, \cite{okte}).

The first step in achieving this goal is to give the {\it complex
   analytic} version of the assignment $[x]\mapsto [x_0]$ explained
above and to prove the analogous remark in this framework.
Therefore, the main object of this article is a holomorphic action
$\alpha:G\times F\ra F$ of a complex reductive group on a complex
manifold $F$. Since we are especially interested in the linear
case (and later in the infinite dimensional case), we will \ub{not}
assume that $F$ is compact. We realized
that  extending the theory of optimal destabilizing OPS to this
situation raises substantial technical difficulties.

First of all, in order to have a good stability condition for a
holomorphic action $\alpha:G\times F\ra F$ one has to fix a K\"ahler
metric $h$  of $F$ which is invariant under a maximal compact
subgroup
$K$ of $G$ and a moment map for the induced $K$-action. Such a data
system $(K,h,\mu)$ provides a generalized maximal weight function
$\lambda:i\kg\times F\ra \R\cup\{\infty\}$.  It is well known (see for instance
Mundet i Riera \cite{mu}) that the stability condition with respect to
$(K,h,\mu)$ can be expressed in terms of the maximal weight function as in the
projective algebraic case. {\it But there is no way to extend this result for
the semistability condition} in the general non-algebraic non-compact case.
Moreover, in the algebraic theory of optimal destabilizing OPSs it is very
important to have a $G$-equivariant maximal weight function, whereas the
choice of a triple $(K,h,\mu)$  only provides a $K$-equivariant one.

In order to solve these difficulties one has to impose a certain
completeness condition on the triple $(K,h,\mu)$, namely {\it
energy completeness} which was introduced in \cite{t} and used in
\cite{lute}. This condition is always satisfied in both compact
and linear case \cite{t}  and also for certain moduli problems on
curves \cite{brte2}.

Moreover, in order to get a $G$-equivariant maximal weight
function $\lambda$, it is convenient to work with an equivalence
class of triples $(K,h,\mu)$ and to show that $\lambda$ extends to
the union of all subspaces of the form $i\kg$. The equivalence is
defined by the natural $G$-action on the set of such triples.
Such an equivalence class will be called a {\it symplectization}
of the action $\alpha$, and it plays the same role as a
linearization of an action in an ample line bundle, in classical
GIT.

An important tool in our proofs will be the linearization theorems of
Heinzner-Huckleberry for Hamiltonian actions \cite{hehu}.

The contents of this article is the following: First we explain
the properties of the maximal weight function associated with an
energy-complete symplectization. Next we prove one of our main
results: the existence and the unicity (up to equivalence) of an
optimal destabilizing element $\xi$ in the Lie algebra of $G$ for
any non-semistable point $f\in F$.    Following the principle
explained in the algebraic case, we next show that the path
$e^{t\xi} f$ converges to a point $f_0$ which is semistable with
respect to a natural action of the  reductive centralizator
$Z(\xi)$ on a certain submanifold of $F$. Fixing the conjugacy
class of $\xi$ one gets a $G$-invariant subset of $F$. The subsets
of this type give a $G$-invariant stratification of $F$, which,
for a large class of actions, is locally finite with locally
Zariski closed strata.   This stratification is the analogue of
the Shatz stratification in the theory of holomorphic vector
bundles (see \cite{ki} for the projective case).

At the end, we
  study  the    optimal destabilizing vectors  of the non-semi\-stable 
objects of two  important  gauge theoretical
moduli problems: holomorphic bundles and holomorphic pairs (bundles
coupled with morphisms with fixed source). Detailed proofs of these results
can be found in \cite{Br}. In this way, we   illustrate our general principle:
\\

{\it In order to get the analogue of the Harder-Narasimhan theorem
for a complex geometric moduli problem one has to give a gauge
theoretical formulation of the problem and to study the optimal
destabilizing vectors of the non-semistable objects.}
\\

This suggests that this principle also holds in the infinite 
dimensional gauge theoretical framework. Details on the gauge 
theoretical
examples and generalizations will appear in a future article.

\section{Background}

\subsection{Symplectization of a holomorphic action and the maximal
weight map $\lambda $}

{\ }\vspace{4mm}

Let us recall some   definitions introduced   in \cite {t}.

Let $G$ be a complex reductive group. One can identify the set
$\Hom(\C^*,G)$ of one
parameter subgroups of $G$ with a subset  $H_{\rm alg}(G)$ of the Lie algebra $\g$
using the map
$$\lambda\mapsto  d_1(\lambda)(1)=\frac{d}{dt}\big|_{t=0}(\lambda(e^t))\ ,$$
where $\lambda$ was regarded here as a map between real manifolds and the Lie algebra of $S^1\subset\C$ was identified with $i\R$.
In  complex non-algebraic
geometry one has to consider  a larger subset of $\g$, and to define
a ''generalized maximal weight function" on this larger set. In the 
non-algebraic non-compact case the optimal destabilizing vectors do 
{\it not} belong
  in general to $H_{\rm alg}(G)$, and this phenomenon occurs even in the simple case of
linear actions (see section \ref{section6}).

\begin{definition}
Let $G$ be a complex reductive group and $\g$ its Lie algebra. We
denote by $H(G)$ the subset of $\g$ consisting of elements $s\in \g$
{\it of Hermitian type}, i.e. of elements which satisfy one of the
following equivalent properties:
\begin{enumerate}
\item There exists a compact subgroup $K\subset G$ such that $s\in i\k$.
\item For every embedding $\rho :G\mapsto GL(r,\C)$ the matrix
$\rho_\star(s)$ is diagonalizable and has real eigenvalues.
\item The closure of the real one parameter subgroup of $G$ defined
by $is\in\g$ is compact.
\end{enumerate}

\end{definition}

This subset is invariant under the adjoint action of $G$ on $\g$; in
general it is not closed. Although $H(G)$ is a subset of $\g$, it {\it cannot}  be
defined intrinsically in terms of the Lie algebra $\g$.

\

One can associate to every $s\in H(G)$ a parabolic subgroup
$G(s)\subset G$ in the following way :
$$G(s):=\{g\in G | \lim_{t\to +\infty} e^{st}ge^{-st} \text{ exists in } G\}.$$
Then $G(s)$ decomposes as a semi-direct product $G(s)=Z(s) \ltimes
U(s)$, where $Z(s)$ is the centralizer of $s$ in $G$ and $U(s)$ is
the unipotent subgroup defined by :
$$U(s):= \{g\in G | \lim_{t\to +\infty} e^{st}ge^{-st}=e\}.$$
We will denote by $\g(s),\mathfrak z(s)$ and $\mathfrak u(s)$ the
corresponding Lie algebras, and by $p_{{\mathfrak z(s)}}$, $p_{{\mathfrak
u(s)}}$ the corresponding projections.

\

Recall the following facts from   \cite {t} :

\begin{proposition}
$\ $

\begin{enumerate}
\item\label{con1} Let $\sigma,s \in H(G)$. The following properties
are equivalent :
\begin{enumerate}
\item $s$ and $\sigma $ are conjugated under the adjoint action of $U(s)$;
\item $s$ and $\sigma $ are conjugated under the adjoint action of $G(s)$;
\item $\sigma \in \g(s)$ and $p_{\mathfrak z(s)}(\sigma )=s$.
\end{enumerate}
\item If one of these conditions is satisfied then $G(s)=G(\sigma
)$. \item The condition in \ref{con1}) defines an equivalence
relation $\sim$ on $H(G)$. \item Let $K$ be a maximal compact
subgroup of $G$. Then $i\mathfrak k\subset H(G)$ is a complete
system of representatives for $\sim$. Mapping $s$ to the
representative in $i\mathfrak k$ of its equivalence class gives a
continuous retraction $\sigma _K: H(G) \to i\mathfrak k$.

\end{enumerate}

\end{proposition}

\begin{example}
Assume that $G=GL(r,\C)$. The data of an equivalence class of $H(G)$
is the data of a pair $(\mathcal F,\lambda)$ where $\mathcal F$ is a
filtration
$$\mathcal F : \{0\} \subset V_1 \subset \cdots \subset V_k=\C^r$$
   and $\lambda$ is an increasing sequence $\lambda_1 < \cdots
<\lambda_k$ of real numbers. An element $s\in \g l(r,\C)$ belongs to
the equivalence class corresponding to $(\mathcal F,\lambda)$ if it
is diagonalisable with spectrum $(\lambda_1,  \cdots ,\lambda_k)$ and
   $$V_i=\bigoplus_{\scriptstyle j=1}^i V_{\lambda_i}\ ,
   $$
   where $V_{\lambda_i}$ is the $\lambda_i$ eigenspace of $s$.   Here,
the parabolic subgroup $G(s)$   is  the subgoup of  matrices
stabilizing the filtration $\mathcal F$.
\end{example}

Following   \cite{t}, we introduce the notion of
\emph{symplectization} of an holomorphic action. A symplectization
of a holomorphic action $\alpha$  plays the same role as  a
linearization of an algebraic action in an ample line bundle in
the classical GIT. This notion will allow us to define a
$G$-equivariant maximal weight function on the set $H(G)$.

\begin{definition}[\cite{lute}]
A {\emph{symplectization}} of the action $\alpha $ is an equivalence
class of triples $(K,h,\mu)$, where $K$ is a maximal compact subgroup
of $G$, $h$ is a $K$-invariant K\"ahler metric on $F$
and $\mu :F \to \mathfrak k\check{ }$ is a moment map for the
$K$-action with respect to the symplectic structure $\omega _h$
defined by $h$.

Two 3-tuples $(K,h,\mu)$ and $(K',h',\mu')$ will be considered
equivalent if there exists $\gamma \in G$ such that :
$$ K'=\Ad_\gamma (K), \ \ h'=(\gamma ^{-1})^\star h, \ \
\mu'=\ad^t_{\gamma ^{-1}}\circ \mu \circ \gamma ^{-1}$$
    \end{definition}

A symplectization of a holomorphic action $\alpha$ should be
regarded as a {\it complex geometric} datum,  which allows one to
define a stability condition independently of the choice of a
maximal compact subgroup of $G$. A triple $(K,h,\mu)\in\sigma$
should be regarded as a {\it symplectic geometric} parameter
compatible with the complex geometric data $\sigma$.
\\

Let $f\in F$ and $u\in \g$. We denote by $c_f^u$ the path in $F$ defined by
$$c_f^u:[0,\infty) \to F\ ,\  c_f^u(t):=e^{tu}f\ .
$$

In order to define the ``maximal weight'' map $\lambda$ associated
with a symplectization, let us introduce the following definition:
\begin{definition}[\cite{lute}]\label{enco}
A symplectization $\sigma $ of the action $\alpha $ will be called
{\emph{energy-complete}} if, for a representative $(K,h,\mu )\in
\sigma$ (and hence for any representative) the following implication
holds: if  $s\in i\mathfrak k$, $f\in F$ and the energy $E_h(c_f^s)$
with respect to the Riemannian metric $h$ is finite,  then $c_f^s$
has a limit as $t\to +\infty$.
\end{definition}

Let $\alpha:G\times V\to V$ be a linear action. A symplectization
of $\alpha$ given by a triple $(K,h,\mu)$, where $h$ is a
Hermitian structure on the vector space $V$,  will be called a
{\it linear symplectization} of $\alpha$.

\begin{remark}
Any linear symplectization and any symplectization of an action on
a compact complex manifold  is  energy-complete \cite{t}. The natural
symplectisation of the action of the complex gauge group on the
configuration space associated to certain moduli problems on curves
is also energy complete \cite{brte2}.
\end{remark}

If we choose a representative $(K,h,\mu)\in \sigma $, we can
associate to every pair $(s,t)\in i\mathfrak k\times \R$ the map
\begin{align*}
   \lambda _t^s : F &\to \R \\
f & \mapsto  \mu ^{-is}(e^{ts}f)
\end{align*}
where we use the notation $\mu^s:=\langle\mu ,s\rangle :F\to \R$ for
any $s\in \mathfrak k$.

It is easy to see that the map $t\mapsto \lambda_t^s(f)$ is
increasing so that one can put
$$\lambda^s(f):= \lim_{t\to +\infty} \lambda_t^s(f) \in \R \cup \{\infty\}.$$

The energy-completeness condition allows one to prove the following
technical result.
\begin{proposition}[\cite{t}]
Assume that $\sigma $ is energy-complete and let $s\in H(G)$. The
map $\lambda^s:F\to \R\cup\{\infty\}$ does not depend on the choice of a
representative $(K,h,\mu)\in \sigma $ with $s\in i\mathfrak k$ and
gives rise to a well defined map
\begin{align*}
\lambda : H(G)\times F &\to  \R\cup\{\infty\} \\
(s,f) &\mapsto  \lambda (s,f)=\lambda^s(f)
\end{align*}
\end{proposition}

The following properties of the map $\lambda $ will be useful in  our study :

\begin{proposition}[\cite{t}]\label{lampro}
Assume that $\sigma $ is energy-complete. The map $\lambda $
introduced above has the following properties :
\begin{enumerate}
\item homogeneity : $\lambda (ts,f)=t\lambda (s,f)$ for any $t\in
\R^+$; \item $\lambda $ is $G$-equivariant: $\lambda (s,f)=\lambda
(\ad_{\gamma^{-1}}(s),\gamma^{-1}.f)$; \item $\lambda $ is $\sim$
invariant : $\lambda ^s(f)=\lambda ^\sigma (f)$ if $s\sim \sigma
$; \item semi-continuity~: \\if $(f_n,s_n)_n \to (f,s)$, then
$\lambda^s(f) \leqslant \liminf_{n\to \infty}\lambda ^{s_n}(f_n)$.
\end{enumerate}
\end{proposition}
\begin{remark} One can work with a similar equivalence relation $\simeq$ on
$H(G)$, using the parabolic subgroups
$$G^-(s):=\{g\in G|\ \lim_{t\to-\infty} e^{ts} g e^{-ts}\ {\rm exists\ in\ }
G\}\ .
$$
It holds
$$\sigma\sim s\Leftrightarrow -\sigma\simeq -s\ .
$$
In general, our map $\lambda$ will {\it not} be $\simeq$ invariant
(see Property \ref{lampro} above).  On the other hand, the "oposite" maximal
weight map given by $\lambda^s_- =\lambda^{-s}$ will be invariant with respect
to this relation. Using $\simeq$ and $\lambda_-$, one will get a completely
parallel theory.

\end{remark}

\subsection{Analytic and symplectic stability}

   {\ }\vspace{4mm}

Let also $\alpha $ be an action of a reductive group $G$ on a complex
K\"ahler manifold $F$, let us choose an  energy-complete
symplectization $\sigma $,  and let $\lambda :H(G)\times F \to \R\cup
\{\infty\}$ be the associated maximal weight map.

We will denote by   $s^\sharp$   the vector field
on $F$ defined by $s$. We will denote  by $\g_f$ the Lie algebra
of the stabilizer of a  point $f\in F$, hence the Lie subalgebra
of $\g$ consisting of those elements $s$ such that $s^\#_f=0$.

\begin{definition} A point $f\in F$ will be called
\begin{enumerate}
\item analytically $\sigma $-semistable if $\lambda^s(f)\geqslant 0$
for all $s\in H(G)$.
\item analytically $\sigma $-stable if it is semistable and $\lambda
^s(f)>0$ for any $s\in H(G)\backslash\{0\}$.
\item analytically $\sigma $-polystable if it is semistable, $\g_f$
is a reductive subalgebra and $\lambda ^s(f)>0$ for every $s$ which
is not equivalent to an element of $\g_f$.
\end{enumerate}
\end{definition}

In this definition we used the following convention: A subalgebra
$\g'\subset \g$ is called a reductive subalgebra if it has the
form $\g'=\mathfrak {k '}^\C$, where $\mathfrak k'$ is the Lie
algebra of a compact subgroup of $G$.  This is more restrictive
than the condition that $\g'$ is isomorphic to the Lie  algebra of
a reductive group.

\begin{remark}
   \item The property of stability (semistability and polystability)
for $f\in F$ depends only on the complex orbit $Gf$ of $f$.
\end{remark}

Note that the proof of this fact for semistable points requires
energy-completeness.

Let us remind the classical definition of (semi)stability for
symplectic actions (see \cite{hi}, \cite{ki}, \cite{helo},
\cite{hehu}). The polystability condition was first introduced in
\cite{ost} in the algebraic framework, as a natural generalization
of the polystability condition for bundles.

\begin{definition}
Let $\sigma $ be a symplectization of the action $\alpha :G\times
F\to F$. A point $f\in F$ is called
\begin{enumerate}
\item symplectically $\sigma$-semistable if, choosing any
representative $(K,h,\mu)\in \sigma $, one has $\overline{G.f}\cap
\mu^{-1}(0)\not= \emptyset$. \item symplectically $\sigma $-stable
if $G.f\cap\mu^{-1}(0)\not=\emptyset$ and $\g_f=\{0\}$. \item
symplectically $\sigma $-polystable if $G.f\cap \mu^{-1}(0)\not=
\emptyset$.
\end{enumerate}
\end{definition}

These conditions do not depend on the chosen representative
$(K,h,\mu)\in \sigma$ and they are obviously  $G$-invariant
conditions wtih respect to $f$.  Note also that the polystability
condition is not open in general.

The following result of Heinzner and Loose (see \cite{helo},
\cite{hehu})  show that one can  always construct a good quotient of the
semistable locus. No  condition on the symplectization is needed.

\begin{theorem}
The set $F^{ss}(\sigma )$ of symplectically $\sigma$-semistable
points is open. Moreover, there is a categorical  quotient
$$F^{ss}(\sigma ) \to Q_\sigma$$
where $Q_\sigma$ is a Hausdorff space with the property that two
$G$-orbits have the same image in $Q_\sigma$ if and only if their
closure contains a common symplectically $\sigma $-polystable orbit.

Choose a representative $(K,h,\mu)\in \sigma $, then every $\sigma
$-polystable orbit intersects $\mu ^{-1}(0)$ along a $K$-orbit and
the induced map
$$\mu ^{-1}(0)/K \to Q_\sigma $$
is a homeomorphism.
\end{theorem}

The following fundamental result links these two notions of stability :

\begin{theorem}[\cite{mu}, \cite{t}]
Assume that $\sigma $ is energy-complete. A point $f$ is
symplectically $\sigma $-stable (polystable) if and only if  it is
analytically $\sigma $-stable (polystable).
\end{theorem}

Our goal here is merely to study the behavior of non semistable
points. So from our point of view, the most important fact is that
the concepts of analytic semistability and symplectic
semistability coincide.  This is a rather difficult technical
result (\cite{t} for details). The main tool  is
the so-called integral of the moment map, whose existence is
assured by the following

\begin{lemma}\label{psi}
Let $(K,h,\mu)$ be a representative of the symplectization $\sigma
$, then there exists a unique smooth function $\Psi : F\times G
\to \R$ with the following properties:
\begin{itemize}
\item $\frac{d}{dt}\Psi (f,e^{ts})=\lambda _t^s(f)$.
\item $\Psi(f,k)=0$ for all $k\in K$.
\item $\Psi (f,gh)=\Psi (f,h)+\Psi (hf,g)$, for all $h,g\in G, f\in F$.
\end{itemize}
\end{lemma}

\begin{proof}
This is a   well-known result (see for instance \cite{mu}).
\end{proof}

\begin{remark}\hfill{\break}
\vspace{-4mm}
\begin{enumerate}
\item The map $t\mapsto \Psi (f,e^{ts})$ is convex for all $s\in
i\mathfrak k, \ f\in F$.

\item  The two following  properties are equivalent :
\begin{itemize}
\item $g\in G$ is a critical point of the map $g\mapsto \Psi (f,g)$;
\item $\mu(gf)=0$.
\end{itemize}
\end{enumerate}
\end{remark}

\begin{theorem}[see  \cite{t}]
Let $(F,h)$ be a K\"ahler manifold, $\alpha : G \times F \to F$ a
complex reductive Lie group action and let $\sigma $ be an energy
complete symplectization for this action.   Then, for any point $f\in
F$ the following properties are equivalent~:
\begin{enumerate}
\item the point $f$ is analytically $\sigma $-semistable; \item
the map $g\to \Psi(f,g)$ associated to any representative
$(K,h,\mu)$ of $\sigma $ is bounded from below over $G$; \item the
point $f$ is symplectically $\sigma $-semistable.
\end{enumerate}

\end{theorem}

Note that energy completeness plays an essential role in the proof.
In the sequel we will speak of $\sigma $-semistability (stability,
polystability) without precising if  the analytical  or symplectical
condition is meant.

\section{The reductive quotient associated to a class of Hermitian
type elements and its canonical action}

We have seen that any equivalence class of elements of Hermitian
type defines a parabolic subgroup $G(\S)$ of $G$.  In this
section, our purpose is to associate to any non trivial
equivalence class $\S$ of Hermitian type elements a new
factorization problem with symmetry group $G(\S)/U(\S)$, where
$U(\S)$ is the unipotent subgroup associated with $\S$. This
quotient is a reductive group.
The new manifold, that we introduce is isomorphic to a submanifold
of $F$, but the identification is not canonical. Then we will show
that for any choice of a symplectization $\sigma$ for the
factorization problem $(F,G,\alpha)$, we may define a natural
symplectization for our new  problem associated to the class $\S$.

\subsection{Natural  action of the canonical reductive quotient}\label{reduquo}

{\ }\vspace{4mm}

   First of all, let us remind that for any $s,s'\in
H(G)$ such that $s\sim s'$, we have $G(s)=G(s')$ and $U(s)=U(s')$
(because $U(s)$ is a normal subgroup in $G(s)$), so that we may
associate to any equivalence class $\S\in H(G)/\sim$ of Hermitian
type elements a unique parabolic subgroup $G(\S)$ of $G$ and a
unique  unipotent subgroup $U(\S)\subset G(\S)$.

For every $s\in H(G)$, let us denote by $$V_s :=\{f\in F \, | \,
(s^\sharp)_f=0\}$$ the zero locus of the vector field $s^\sharp$.
Locally this set consists of fixed points under the action of the
compact torus $T=\overline{\{e^{its}\, | \, t\in \R\}}$ and
therefore, using the slice theorem \cite{orl}, we see that $V_s$
is a smooth submanifold of $F$, in general {\it not} of pure dimension. Being
the vanishing locus of  the holomorphic tangent field associated
with $s^\sharp$, it  inherits   a structure of complex manifold
(of possibly non-pure dimension).

Let us  now remark that for any $s$ and $s'$ in $H(G)$, if $s\sim
s'$ then there exists $u\in U(s)=U(s')$ such that $s'=\ad_u(s)$
and we have an associated isomorphism $\alpha (u): V_s
\xrightarrow{\simeq} V_{s'}=V_{\ad_u(s)}$.

One can easily prove that the element $u\in U(s)=U(s')$ such that
$s'=\ad_u(s)$ is unique, so that one gets a canonical
identification $V_s\simeq V_{s'}$. Indeed, if there exists two
elements $u,v\in U(s)$ such that $\ad_u(s)=\ad_v(s)=s'$, we get
$w=v^{-1}u \in U(s)$ and $\ad_w(s)=s$. Then we have, $w\in
Z(s)\cap U(s) = \{e\}$ so that the induced isomorphism is the
identity.

Therefore we can associate to any non trivial   equivalence class
$\mathcal S$ of $H(G)$ a canonically defined complex manifold
$$ \V (\S) := \{ \coprod_{s \in \S} V_s\} / \sim $$
where $\sim$ is induced by the previous identifications. One has, for every
$s\in{\cal S}$, a natural identification $\V (\S)\simeq V_s$.
\begin{remark} One can wonder why is it important to consider the copy $\V
(\S)$ of $V_s$. The reason is the following: we will associate to any
non-semistable point $f\in F$   a  well defined   {\it class} $\S_f\in
H(G)/\sim$   of so-called optimal destabilizing vectors, but there is {\it no
way} to associate a well defined such destabilizing vector.
Therefore, we get a well defined assignment    $f\mapsto  \V(\S_f)$.  In the
next step we will see that $f$ also defines canonically a point $f_0\in
\V(\S_f)$ which is semistable with respect to a certain (again
canonically associated) symplectization.
\end{remark}

\

The action $\alpha $ induces an action of the parabolic group $G(\S)$
over the complex manifold $\V(\S)$ defined by :
\begin{align*}
G(\S)\times \V(\S) & \to \V(\S) \\
(g, [x]) & \mapsto [g(x)]
\end{align*}
where $x\in V_s$ and $g(x)\in V_{\ad_g(s)}$ for any $s\in \S$.

Of course $G(\S)$ is not reductive but it is easy to see, using the
definition of $\V(\S)$,  that  the unipotent  subgroup $U(\S)$ acts
trivially on $\V(\S)$. So that we get a well-defined action
$$\alpha_\S : G(\S)/\U(\S) \times \V(\S) \to \V(\S)$$
of the canonical reductive quotient $G(\S)/U(\S)$.

Let us remark that if we choose any representative $s\in \S$, the
action of $G(\S)/U(\S)$ over the representative $V_s$ of $\V(\S)$
is just the induced action of the reductive Lie group $Z(s)$ over
the complex submanifold $V_s \subset F$.

\

\subsection{A natural symplectization for the action
$\alpha_\S$}\label{natsymp}

{\ }\vspace{4mm}

An \emph{$\ad$-invariant inner product of Euclidian type} on the
Lie algebra $\mathfrak{g}$ is an $\ad_G$-invariant non-degenerate
complex symmetric   bilinear form $h$ on $\mathfrak g$ which
restricts to an inner product on a subspace of the form
$i\mathfrak k$ (and hence on any subspace of this form as any two
such subspaces are  conjugated).

The data of such   an inner product is equivalent to the data of:\\
--  a multiple  of the Killing form $k_{\mathfrak{s}}$ of each simple
summand $\mathfrak{s}$ of the semisimple part $\mathfrak g^s$ of
$\mathfrak g$ and  \\
-- an inner  product on $i\mathfrak t_0$, where $\mathfrak t_0$ is
the Lie algebra   of the maximal compact subgroup of the complex
torus $\mathfrak z(\mathfrak g)$.\\

We fix such an inner product $\langle\cdot,\cdot\rangle$ on our Lie
algebra $\mathfrak{g}$.\\

For any choice of a symplectization $\sigma $ of the
factorization problem $(F,G,\alpha)$, we may define   a canonical
symplectization for the action $\alpha_\S$ in the following way:
let  $\rho =(K ,h ,\mu)\in \sigma $, and  let us  take the unique
representative $s\in i\k \cap \S$ and the corresponding copy
$V_s$ of $\V(\S)$, then we can define an associated
symplectization of $V_s\simeq \V(\S)$ using the triple $$\rho
_\S:=(K \cap Z(s), h|_{V_s}, i^\star(\mu |_{V_s}) +\tau )$$ where
$i: \k \cap \z(s) \hookrightarrow \k$ is the inclusion and $\tau $
is the locally constant ${\z}(s)^\vee$-valued function over $V_s$
defined by
$$\tau (x)= - (\mu^{-is}(x)) \langle is,\cdot\rangle\ .$$

To see that the map  $\tau$ above is indeed locally constant, note
that the map $x \to
\mu^{-is}(x)$ is locally constant over $V_s$ since
$$d\mu^{-is}(\cdot)=\omega_h((-is)^\sharp,\cdot)=h(s^\sharp,\cdot)=0\ .$$

   The reason for this particular choice of the moment map
will appear later  in  section \ref{stratquo}.

This definition is coherent with the identifications defined above
: let $\rho'=(K',h',\mu')\in \sigma $ be another representative of
$\sigma $ and let $s'\in i\k' \cap \S$.  Then,  there exists $u\in
U(\S)$ such that $\ad_u(s')=s$. The application $\alpha (u)$
defines an isomorphism from $V_{s'}$ onto $V_s=V_{\ad_u(s)}$ and
conjugates $\rho'$ to another representative
$\rho''=u_\star(\rho')\in \sigma $ defined by $$\rho''=
(K'',h'',\mu'')=  (\Ad_u(K'), (u^{-1})^\star (h'), \ad_{u^{-1}}^t
\circ \mu' \circ u^{-1}).$$ It is sufficient to show that $\rho
_\S$ and ${\rho''}_\S$ define the same symplectization for the
action $\alpha_\S$. Let us remark that, by the definition of
$\tau$, ${\rho''}_\S=u_\star(\rho'_\S)$. Moreover,  we have $s\in
i\k\cap i\k''$ and there exists $\gamma \in G$ which conjugates
$\rho $ and $\rho''$, i.e. $$\Ad_\gamma (K'')=K, \
h=(\gamma^{-1})^\star(h''), \ \mu=\ad^t_{\gamma ^{-1}}\circ \mu''
\circ \gamma^{-1}.$$ We now use the following lemma :

\begin{lemma}
Let $K$ be a maximal compact subgroup of $G$ and let $g\in G$,
$s\in \k$ such that $\ad_g(s)\in \k$. Then the decomposition
$g=kl$, where $l\in exp(i\k)$, $k\in K$ satisfies $\ad_l(s)=s$,
i.e. $l\in Z(s)$.
\end{lemma}

\begin{proof}
Let us decompose $g$ as $g=kl$ with $l\in exp(i\k)$ and $k\in K$. Then
$$ \gamma :=\ad_l(s)=\ad_{k^{-1}}(\ad_g(s)) \in \k.$$
If we choose an embedding $G\hookrightarrow GL(r,\C)$ mapping $K$ to
$U(r)$, then the image of $l$ is Hermitian with positive eigenvalues,
whereas the images of $s$ and $\gamma $ are anti-Hermitian. We have :
$$ \ad_l(s)^\star=-\ad_{l^{-1}}(s)=\gamma^\star = -\gamma =-\ad_l(s),$$
hence $\ad_{l^2}(s)=s$. This implies that the eigenspaces of $l^2$
and hence of $l$ are invariant under $s$, so that one also has
$\ad_l(s)=s$.
\end{proof}

Therefore, since $s\in i\k$ and $\ad_\gamma(s)\in i\k$,  we have
the decomposition $\gamma =kl$, $k\in K$ and $l\in Z(s)$ so that
$$K''=\Ad_{l^{-1}k^{-1}}(K)= \Ad_{l^{-1}}(K),$$
$$h''=\gamma^\star(h)=l^\star(h)$$
   because  $h$ is by definition $K$-invariant and
$$\mu''=\ad_l^t\circ \ad_k^t\circ \mu\circ k \circ l=\ad_l^t\circ \mu
\circ l$$ because a moment map is always $K$-equivariant. We conclude
that $\rho $ and $\rho''$ are conjugated by an element of $l\in
Z(s)$. One has $s''=s$ because they are both representatives in
$i\kg''$ of $\S$,  therefore $\tau$ and $\tau''$ are conjugated, so
that  the two induced triple $\rho_\S$ and $\rho''_\S$ are equivalent
for the action $\alpha_\S$.

In the sequel, we will denote by $\sigma_\S$ this natural
symplectization for the factorization problem
$(\V(\S),G(\S)/U(\S),\alpha_\S)$.

\section{optimal destabilizing vector for a non semistable
point}\label{secopti}

In this section we will associate to every non $\sigma$-semistable
point $f\in F$, an optimal destabilizing element $s\in H(G)$
which minimize the weight  function $\lambda (.,f)$. We will also
see that this element is unique up to equivalence.

\

So, let us consider a holomorphic action $\alpha :G\times F\to F$ of
a reductive group $G$ on the K\"ahler manifold $F$. We choose a
symplectization $\sigma $ for this action and we assume in the sequel
that $\sigma $ is energy-complete (see def. \ref{enco}) so that  the
map $\lambda :H(G) \times F \to \R\cup\{\infty\}$ is well defined.

\

Fix again an \emph{$\ad$-invariant inner product of Euclidian
type} $\langle\cdot,\cdot\rangle$ on $\mathfrak{g}$. Such a
structure gives a well defined real application $\|\cdot\|
:H(G)\to \R $ defined by $\|s\|=\sqrt{\langle
s,s\rangle}$ (in fact all the elements of  $H(G)$  lie
in a  Lie algebra of the form $i\mathfrak k$ for a certain maximal
compact subgroup $K$, on which $\langle\cdot,\cdot\rangle$ is a
scalar product). Let us remark that $\langle\cdot,\cdot\rangle$ is
constant on the equivalence classes of $H(G)$, so that we may
speak of a ``normalized class'' $\S$.

\

We consider in this section a given $\sigma $-non semistable point
$f\in F$ and we set
    $$\lambda _{\rm inf}:= \inf_{\substack{s\in H(G)\\ \|s\|=1}}
\lambda(s,f).$$
Let us remark that this lower bound is not $-\infty$ as
$$\lambda(s,f) \geqslant \lambda _0(s,f) + E_h(c_f^s) \geqslant
\lambda _0(s,f)=\langle\mu(f),-is\rangle$$
Let us define the set of normalized destabilizing elements of $f$
: $$\Lambda_f :=\{\xi \in H(G) \ | \ \|\xi\|=1 \text{
and } \lambda(\xi,f)=\lambda_{\rm inf}\}.$$

\begin{theorem}\label{maxdesele}
Let $f\in F$ be a non $\sigma$-semistable point. Then $\Lambda_f$ is
non empty and consists of exactly   a normalized  equivalence class
$\S_f\subset H(G)$.
\end{theorem}


\begin{proof}$ $
\begin{lemma}[Existence]\label{conju}$ $

\begin{enumerate}
\item If $s\in \Lambda_f$ and
$s'\in H(G) $ with $s'\sim s$ then $s'\in \Lambda _f$.
\item $\Lambda_f\not= \emptyset$.
\end{enumerate}
\end{lemma}

\begin{proof}
The first point follows directly from the equivariance properties
of $\lambda$ (see prop. \ref{lampro}) and the $\ad$-invariance of
$\langle\cdot,\cdot\rangle$.

For the second point, let us fix a maximal compact subgroup $K$ of
$G$. Then we know that $i\mathfrak k \subset \mathfrak g$ is a
complete system of representatives for $\sim$.  By invariance, the
application $\lambda$ restricts to  a map
$\tilde{\lambda}:i\mathfrak k \to \R\cup\{\infty\}$. Take now a sequence
$(s_n)_n\in H(G)$ such that $\lambda (s_n,f)$ converges
to $\lambda _{\rm inf}$ and $\|s_n\|=1$ for all $n$. We take $\tilde
s_n$ to be the representative in $i \mathfrak k$ which is in the
same equivalence class as $s_n$. We still have $\|\tilde s_n\|=1$
and $\tilde \lambda (\tilde s_n)\to \lambda _{\rm inf}$. Now
$i\mathfrak k$ is a closed finite dimensional vector space in
$\mathfrak g$ so that its unit sphere is compact. Thus, we can
extract a converging subsequence $\tilde s_m \to \tilde s$. Now,
the semi-continuity property of $\lambda $ (prop. \ref{lampro})
implies $$\tilde \lambda (\tilde s,f)\leqslant \liminf_{n\to
\infty}\lambda (s_m,f) = \lambda_{\rm inf},$$ i.e all the elements of
the class $\tilde s$ are elements of $\Lambda_f$.
\end{proof}

\begin{lemma}[Unicity]\label{unicity}
The optimal destabilizing element is unique up to equivalence :
$$ \exists \, \xi \in H(G) \text{ s.t. } \Lambda _f=\{s\in H(G) \ | \
\xi \sim s\}= \S(\xi).$$
\end{lemma}

\begin{proof}
Let us choose a representative $(K,h,\mu)\in \sigma $ and let
$\Psi : F\times G\to \R$ the associated integral of the moment map
(see prop. \ref{psi}). We must prove that there exists only one
maximal element in $i\mathfrak k$.
$i\mathfrak k \cap \Lambda _f$.

Our first step is to prove the result when $K=T$ is a real torus.
\begin{lemma}\label{unitorlin}
If $K=T$ is a real torus then there exists
a unique $\xi_T(f) \in i\mathfrak t$ such that $\Lambda_f\cap
i\mathfrak t=\{\xi_T(f)\}$.
\end{lemma}

\begin{proof}
The proof is based on the following lemma :
\begin{lemma}
The map $\Phi_f : i\mathfrak t \to \R$ defined by $\Phi _f (s)=\Psi
(f,e^s)$ is convex on $i\mathfrak t$.
\end{lemma}

\begin{proof}
This is a well-known property of $\Psi $ that the maps $t\mapsto \Psi
(f,e^{ts})$ are convex for all $s\in i\mathfrak t$ (prop \ref{psi}).
Let $\xi, s \in i\mathfrak t$, then using  the fact that $\xi $ and
$s$ commute, we get :
$$\Psi (f,e^{\xi +ts})=\Psi (f,e^{\xi }e^{ts}) = \Psi (f,e^\xi )+
\Psi (e^\xi  f,e^{ts})$$
so that $t \mapsto \Phi (\xi +ts)$ is convex for every $s\in
i\mathfrak t$. To conclude we use the following easy lemma:
\begin{lemma}
Let $f:U\subset \R^n \to \R$ a smooth function such that for all
$x_0,x\in \R^n$ the map $t\mapsto f(x_0+tx)$ is convex. Then $f$ is
convex on $\R^n$.
\end{lemma}
\end{proof}

By definition we have
$$\lambda(\xi,f)=\lim_{t\to +\infty} \frac{d}{dt}\phi_f(t\xi)$$
Assume that  there exist two distinct optimal destabilizing
elements $\xi_1, \xi_2 \in i\mathfrak t\cap \Lambda_f$ and let
$\xi=\frac{\xi_1+\xi_2}{2} \in i\mathfrak t$. Of course we have
$\|\xi\|<1$.  The convexity of the function $\phi_f$ implies  that
$\phi_f(t\xi)\leqslant \frac{1}{2}(\phi_f(t\xi_1)+\phi_f(t\xi_2))$
for all $t\in \R$. We get
\begin{align*}
\frac{\phi_f(\theta \xi) -\phi_f(t\xi)}{\theta - t} \leqslant &
\frac{\frac{1}{2}(\phi_f(\theta \xi_1)+\phi_f(\theta \xi_2)) -
\frac{1}{2}(\phi_f(t \xi_1)+\phi_f(t \xi_2))}{\theta-t}  + \\
   \ & \ \  \ \ \ \ \ \ \frac{\frac{1}{2}(\phi_f(t \xi_1)+\phi_f(t \xi_2)) -
\phi_f(t\xi)}{\theta-t}
\end{align*}
and so
$$\hspace*{-.5cm}\limsup_{\theta \to +\infty} \bigg[
\frac{\phi_f(\theta \xi) -\phi_f(t\xi)}{\theta -
t}\bigg] \leqslant\limsup_{\theta \to \infty} \Bigg[
\frac{\frac{1}{2}(\phi_f(\theta \xi_1)-\phi_f(t \xi_1)) +
\frac{1}{2}(\phi_f(\theta \xi_2)-\phi_f(t
\xi_2))}{\theta-t}\Bigg]$$ The regularity and the convexity of
$\Psi$ implies that for all $\theta,t \in \R$
$$\frac{d}{ds}_{|_{s=\theta}}\Psi(f,e^{s \xi}) \geqslant
\frac{\phi_f(\theta \xi) -\phi_f(t\xi)}{\theta -
t} \geqslant \frac{d}{ds}_{|_{s=t}}\Psi(f,e^{s \xi})$$ thus we
have
$$ \hspace*{-.6cm}\lim_{t\to +\infty} \frac{d}{dt}\Psi(f,e^{t \xi}) \leqslant
\limsup_{t \to +\infty} \bigg[\limsup_{\theta \to
+\infty}\bigg(\frac{\frac{1}{2}(\phi_f(\theta \xi_1)-\phi_f(t
\xi_1))}{\theta-t}+\frac{\frac{1}{2}(\phi_f(\theta \xi_2)-\phi_f(t
\xi_2))}{\theta-t}\bigg)\bigg]$$
$$= \frac{\lambda_{\rm inf}+\lambda_{\rm inf}}{2}$$
We deduce from this that
$$\lambda(\frac{\xi}{\|\xi\|},f)=
\frac{\lambda(\xi,f)}{\|\xi\|} \leqslant
\frac{\lambda_{\rm inf}}{\|\xi\|} < \lambda_{\rm inf}\ ,$$
because $\|\xi\|<1$. This leads to a
contradiction.
\end{proof}

\begin{remark} Note that in this argument one essentially needs
the fact that $f$ is non \underbar{semi}stable (i.e. $\lambda_{\rm
inf}<0$).
\end{remark}

   Let us now come back to our main proof for an arbitrary compact lie
group   $K$.
\begin{lemma}
Let $f$ a non $\sigma $-semistable point and $\xi \in \Lambda_f$. Let $T$ be
a maximal torus in $G(\xi)$. Then $f$ is   non semistable with
respect to the induced symplectization of the  $T^\C$-action, and
$\xi$ is conjugated to $\xi_T(f)$ by an element of $G(\xi)$.
\end{lemma}

\begin{proof}
Let $S$ be a maximal torus of $G(\xi)$ whose Lie algebra
$i\mathfrak s$ is containing $\xi$. All maximal tori of $G(\xi)$
are conjugated to each over, so there exists $p\in G(\xi)$ such
that $\Ad_p(S)=T$ then we have $\ad_p(\xi) \in i\mathfrak t \cap
\Lambda_f$ (see proposition \ref{conju}). We deduce from this that
$f$ is $T^\C$ non semistable  and $\ad_p(\xi)$ is an optimal
destabilizing element with respect to the $T^\C$-action.
Therefore, from the previous unicity lemma \ref{unitorlin},
$\xi_T(f)=\ad_p(\xi)$.
\end{proof}

Now we make use of the following well-known lemma (see \cite{slo})~:
\begin{lemma}
Let $P$ and $P'$ be parabolic subgroups of $G$. Then there is a
maximal torus $T$ of $G$ contained in the intersection $P\cap P'$.
\end{lemma}

We use now the same method as in the algebraic case \cite{slo}.

\begin{lemma}
Let $f$ be a non semistable point and let $\xi_1$ and $\xi_2$ be
two elements of $\Lambda_f$. Then we have $G(\xi_1)=G(\xi_2)$ and
there exists $p\in G(\xi_1)$ such that $\ad_p(\xi_1)=\xi_2$, i.e.
$\xi_1 \sim \xi_2$.
\end{lemma}

\begin{proof}
Let $T$ be  a maximal torus contained in $G(\xi_1)\cap G(\xi_2)$.
By the previous lemma, there exists $g\in G(\xi_1)$ and $k\in
G(\xi_2)$ such that $\ad_g(\xi_1)=\xi_T(v)=\ad_k(\xi_2)$. Then we
get $G(\xi_T(v))=\Ad_g(G(\xi_1))=G(\xi_1)$ and the same thing for
$G(\xi_2)$. We get $G(\xi_1)=G(\xi_2)$ and
$\xi_2=\ad_{(k^{-1}g)}\xi_1$, so that $\xi_1 \sim \xi_2$.
\end{proof}
This concludes the proof of the proposition.
\end{proof}
\end{proof}

\section{Associating a semistable point to a non-semistable one. The
Shatz stratification associated with a Hamiltonian
action}\label{stratquo}
Using the results of the two previous sections, we will show here
that it is possible to associate naturally to any non $\sigma
$-semistable point a semistable point for the new factorization
problem defined in section \ref{reduquo} and thus a point in the
associated Hamiltonian quotient. This leads to the stratification
of $F$ by $G$ invariant subsets described in the main thm.
\ref{principal}.  The main stratum is the semi-stable locus and it
is open. The other strata are obtained by fixing the conjugacy
class  (with respect to the $\ad_{G}$-action on $H(\mathfrak{g})$)
of the optimal destabilizing element.

\

Let $(F,G,\alpha)$ be a factorization problem with an associated
symplectization $\sigma $. Let us choose a $\ad_G$-invariant inner
product of Euclidian type over $\g$. Using thm \ref{maxdesele}, we
define for any normalized equivalence class $\S \in H(G)/\sim$
the following subset of $F$ :
$$\mathcal Z_\S := \{f\in F \, | \, f \text{ is non $\sigma $-semistable
and } \Lambda_f = \S\}$$
Therefore $\mathcal Z_\S$ is the locus of points with optimal
destabilizing class $\S$.

Now let us choose $f\in\mathcal Z_\S$, and  let us fix any
representative  $\rho=(K,h,\mu)\in \sigma$ and take the unique
representative $\xi^\rho_f \in \Lambda_f\cap i\mathfrak k$ of
$\S$. The point $f$ being non $\sigma $-semistable, we get from
the formula
$$\lambda(\xi^\rho_f,f)=\lambda_0(\xi^\rho_f ,f) + E_h(c_f^{\xi^\rho_f} )$$
that $
E_h(c_f^{\xi^\rho_f} ) < +\infty$. The symplectization $\sigma $ is
supposed to be energy complete, so that there exists a limit
element $f_0^\rho=\lim_{t\to \infty} e^{t\xi^\rho_f }f \in F\cap
\overline{Gf}$.

\

The point $f_0^\rho $ lies  in the vanishing subset $V_{\xi
^\rho_f}=V((\xi ^\rho_f)^\sharp)$. For another choice $\rho'\in
\sigma $ we get  $\xi^{\rho'_f}=\ad_u(\xi^\rho_f)$, where $u\in
U(\S)$, and   $f_0^{\rho'}=\alpha (u)(f_0^\rho) \in
V_{\xi^{\rho'_f}}$. Thus, we obtain a well defined point $f_0\in
\V(\S)$ canonically associated to $f$. Our claim is the following
:

\begin{theorem}\label{harder}
Let $f\in F$ be a  non $\sigma $-semistable point, and $\S=\Lambda_f$
the class of its optimal destabilizing element. Then the canonically
associated point $f_0 \in \V(\S)$ is $\sigma_\S$-semistable for the
action of the canonical reductive quotient $G(\S)/U(\S)$ over
$\V(\S)$.
\end{theorem}

\begin{proof}
For our purpose, we may fix a representative $\rho=(K,h,\mu)\in
\sigma $ and the element $\xi^\rho_f \in \Lambda_f\cap i\k$. Let
us remark first that, by definition,
$$\lambda_{\rm inf}= \lim_{t\to +\infty}
\mu^{-i{\xi^\rho_f}}(e^{t\xi^\rho_f}f)=
\mu^{-i{\xi^\rho_f}}(f_0^\rho).$$
Let  $$\rho_\S= (K\cap Z(\xi^\rho_f), h|_{V_{\xi^\rho_f}}, \mu'=
i^\star (\mu|_{V_{\xi^\rho_f}}) +\tau)$$ be the associated triple
representing $\sigma _\S$ (see \ref{natsymp}), then $\tau$ is
given on the connected component containing $f_0^\rho $ by $$\tau=
- \lambda_{\rm inf}\langle i\xi^\rho_f,\cdot\rangle \ .$$ Let
$\lambda'$ be the map associated to the symplectization
$\sigma_\S$. We must show that $f_0^\rho $ is $\sigma_\S$
semistable.

An element $s\in i\mathfrak k\cap \mathfrak z(\xi^\rho_f)$ has an orthogonal
decomposition   $s=\beta\xi^\rho_f + s^\bot$.

  From now on, we assume that $\lambda'(s,f_0^\rho)<0$ and we will   get a
contradiction.

Let $\xi_\varepsilon= {\xi^\rho_f} + \varepsilon s^\bot$ for $\varepsilon
>0$. Then we get
\begin{align*}
\lambda (\frac{ \xi_\varepsilon}{\|\xi_\varepsilon\|},f)  &=
\frac{\lambda(\xi_\varepsilon,f)}{\|\xi_\varepsilon\|} \\
&= \frac{\lim_{t\to \infty} \mu^{-i\xi_\varepsilon}(e^{t{\xi}
_\varepsilon }f)}{\|{\xi^\rho_f} +\varepsilon s^\bot\|} \\
& = \frac{\lim_{t\to\infty} \mu^{-i{\xi^\rho_f} }(e^{t{\xi}_\varepsilon }f)
+ \varepsilon \lim_{t\to \infty} \mu^{-is^\bot}(e^{t{\xi}_\varepsilon
}f)}{\|{\xi^\rho_f} +\varepsilon s^\bot\|}\ . \end{align*}

   So we
are reduced to study the orbit of $f$ under the flow of
$\xi_\varepsilon^\sharp$. We begin with the remark that the
hypothesis $\lambda'(s,f_0^\rho)<0$ implies that
$E_h(c_{f_0^\rho}^{s^\bot})<\infty$, so that, the action being energy complete,
we know that the curve $c_{f_0^\rho}^{s^\bot}$ converges to some point $f_1
=\lim_{t\to\infty} e^{ts^\bot}f_0^\rho \in F\cap V_{\xi^\rho_f}$.

The main point of the proof is the following \\ \\
{\bf Claim:}
{\it For any sufficiently small $\varepsilon>0$ , the orbit of $f$ under the
one-parameter subgroup generated by
${\xi}_\varepsilon $ converges to $f_1$, i.e. $\lim_{t\to\infty}
e^{t{\xi}_\varepsilon}f = f_1$.}

{\samepage\begin{proof} (of the Claim) We consider first the compact torus
$$T:=\overline{\{e^{it{\xi^\rho_f}}e^{i\theta s^\bot}\, | \, t,\theta
\in \R\}}\subset K
$$
and the induced action $T^\C\times F\rightarrow F$ of its
complexification $T^\C\subset G$.

Now we use a fundamental result Heinzner and Huckleberry, which
allows us to "linearize" this action around $f_1$. Indeed, up to a
modification of the moment map $\mu_T$ by a constant in $\mathfrak
t=\mathfrak z(\mathfrak t)$, we may always assume that
$\mu_T(f_1)=0$. Now following \cite{hehu} (p. 346), we may find an
open $T^\C$-stable Stein neighborhood of $f_1$. Using the fact
that $T^\C$ is reductive, we can apply Theorem 3.3.14 in
\cite{hehu} and get the existence of an open $T^\C$-invariant
Stein neighborhood $U$ of $f_1$, a linear representation
$\rho:T^\C\times V\rightarrow V$ and a closed $T^\C$-equivariant
embedding $a: U\rightarrow V$. Since $U$ is open and
$T^\C$-invariant, it follows easily that it contains the points
$f_0^\rho$ and $f$. Put $v_1:=a(f_1)$, $v_0:=a(f_0^\rho)$,
$v:=a(f)$.

We decompose $V$ as
$$V=\bigoplus_{\chi\in R} V_\chi\  ,$$
where $R\subset{\rm Hom}(T,S^1)$ and $\rho(t)|_{V_\chi}=\chi(t){\rm
id}_{V_\chi}$ for all $t\in T$.

Since $\lim_{t\rightarrow\infty} e^{t{\xi^\rho_f}} f=f_0^\rho$, we deduce  that
   $$v=v_0+v_- $$
   where
   $$v_0\in\bigoplus_{d_e(\chi)({\xi^\rho_f})=0} V_\chi\ ,\
v_-\in\bigoplus_{d_e(\chi)({\xi^\rho_f})<0} V_\chi\  .$$
   For sufficiently small $\varepsilon>0$ we get that
$d_e(\chi)({\xi^\rho_f}+\varepsilon s^\bot)<0$ for all $\chi\in R$
for which $d_e(\chi)({\xi^\rho_f})<0$.

For such $\varepsilon$ we get
    that $$\lim_{t\to \infty}e^{t({\xi^\rho_f}+\varepsilon
s^\bot)}v=\lim_{t\to \infty}e^{t({\xi^\rho_f}+\varepsilon
s^\bot)}v_0 =\lim_{t\to \infty}e^{t\varepsilon s^\bot}v_0=v_1\
.$$\end{proof}}

According to the above claim, if $\varepsilon $ is sufficiently
small, our  computation gives
$$\lambda (\frac{ \xi_\varepsilon}{\|\xi_\varepsilon\|},f) =
\frac{\mu^{-i{\xi^\rho_f} }(f_1) + \varepsilon
\mu^{-is^\bot}(f_1)}{\|{\xi^\rho_f} +\varepsilon s^\bot\|}$$
Using the same methods  as before we have
$$ \mu^{-i{\xi^\rho_f} }(f_1) = \mu^{-i{\xi^\rho_f} }(f_0^\rho) =
\lambda ({\xi^\rho_f} ,f) =\lambda_{\rm inf} $$
and moreover, by Lemma \ref{lll} below, we get
$$ \mu^{-is^\bot}(f_1) = {\mu'}^{-is^\bot}(f_1) = \lambda
'(s^\bot,f_0^\rho) <0$$
We obtain
$$ {\frac{d}{d\varepsilon }}|_{\varepsilon =0}(\lambda (\frac{
\xi_\varepsilon}{\|\xi_\varepsilon\|},f))= \mu^{-is^\bot}(f_1) <0.$$
Thus, by taking $\varepsilon $ small enough, we get a normalized
element $\frac{\xi_\varepsilon }{\|\xi_\varepsilon \|} \in
H(G)$ with $\lambda (\frac{
\xi_\varepsilon}{\|\xi_\varepsilon\|},f)<\lambda_{\rm inf}$ which
is a contradiction.
\end{proof}
\begin{lemma}\label{lll}
If $s=\beta\xi^\rho_f + s^\bot$ then $\lambda'(s,f_0^\rho
)=\lambda'(s^\bot, f_0^\rho )$.
\end{lemma}

\begin{proof}

We have
\begin{align*}
   \lambda'(s,f_0^\rho)= &  \lim_{t \to \infty} {\mu'}^{-is}(e^{ts}f_0^\rho) \\
                  = & \lim_{t \to \infty} \mu^{-is}(e^{ts}f_0^\rho) -
\lambda_{\rm inf}\langle{\xi^\rho_f},s\rangle \\
                  = & \lim_{t \to \infty}
(\mu^{-is^\bot}(e^{ts}f_0^\rho) +
\mu^{-i\beta{\xi^\rho_f}}(e^{ts}f_0^\rho)) -
                  \lambda_{\rm inf}\langle{\xi^\rho_f},s\rangle
\end{align*}
Now keep in mind that ${\xi^\rho_f}$ and $s$ commute so that
$$ e^{st}f_0^\rho = e^{s^\bot t}(e^{\beta {\xi^\rho_f}
t}f_0^\rho)=e^{s^\bot t}f_0^\rho. $$
we get
\begin{align*}
\lambda'(s,f_0^\rho) = &  \lim_{t \to \infty}
\mu^{-is^\bot}(e^{ts^\bot}f_0^\rho) +  \lim_{t \to
\infty}\mu^{-i\beta{\xi^\rho_f}}(e^{ts}f_0^\rho) - \lambda_{\rm
inf}\langle{\xi^\rho_f},s\rangle \\
                  = & \lim_{t \to \infty}
{\mu'}^{-is^\bot}(e^{ts^\bot}f_0^\rho) +  \lim_{t \to
\infty}\mu^{-i\beta{\xi^\rho_f}}(e^{ts}f_0^\rho) - \lambda_{\rm
inf}\beta \\
                  = & \lambda'(s^\bot,f_0^\rho) + \beta\lim_{t \to
\infty}\mu^{-i{\xi^\rho_f}}(e^{ts}f_0^\rho) - \lambda_{\rm inf}\beta
\end{align*}
Note that
$$ \mu^{-i{\xi^\rho_f}} (e^{ts}f_0^\rho) = \mu^{-i{\xi^\rho_f}}
(f_0^\rho) + \int_0^t \frac{d}{d\tau}\mu^{-i{\xi^\rho_f}}(e^{\tau
s}f_0^\rho) d\tau $$
Using the definition of the moment map, we get
$$ \frac{d}{d\tau}\mu^{-i{\xi^\rho_f}}(e^{\tau s}f_0^\rho)=
d(\mu^{-i{\xi^\rho_f}})(\mathfrak v_\tau) = \omega_h
(-i{\xi^\rho_f}^\sharp,\mathfrak v_\tau) =
h({\xi^\rho_f}^\sharp,\mathfrak v_\tau)$$
where $\mathfrak v_\tau$ is the speed vector along the curve
$c_{f_0^\rho}^s$. But the vector field ${\xi^\rho_f}^\sharp$ vanishes
identically along the curve  $c_{f_0^\rho}^s$, because
$$e^{t{\xi^\rho_f}}(e^{\tau s}f_0^\rho) = e^{\tau
s}(e^{t{\xi^\rho_f}}f_0^\rho) =e^{\tau
s}f_0^\rho$$
   so that each point $c_{f_0^\rho}^s(\tau)$ of the curve is a
fixed point of the flow of the vector field ${\xi^\rho_f}^\sharp$.
We get $h ({\xi^\rho_f}^\sharp,\mathfrak v_\tau)=0$ and
$\mu^{-i{\xi^\rho_f}} (e^{ts}f_0^\rho)= \mu^{-i{\xi^\rho_f}}
(f_0^\rho)=\lambda_{\rm inf}$. The above formula shows that
$$\lambda'(s,f_0^\rho)=\lambda'(s^\bot,f_0^\rho).$$
\end{proof}

\begin{corollary}
The subsets $\mathcal Z_\S$ are $G(\S)$-invariant and there is a
natural quotient map $\mathcal Z_\S \to Q_{\sigma_\S}$ where
$Q_{\sigma_\S}$ denotes the Hamiltonian quotient associated to the
factorization problem
   $(\V(\S),G(\S)/U(\S))$ and to the symplectization $\alpha_\S$.
\end{corollary}

\begin{proof}
The invariance is a direct consequence of the $\ad$-invariance
properties of $\lambda$ (see prop. \ref{lampro}).
\end{proof}

To get a $G$ invariant stratification  we have to glue these
subsets together in the following way : $H(G)$ is $\ad_G$
invariant and we denote by $\Sigma_{\ad}$ the set of all orbits
for this action. Then for any non trivial orbit $\delta \in
\Sigma_{\ad} $, we define
$$\mathcal X_\delta :=\{ f \in F \,  |\, f\ {\rm non-semistable},\ \S_f \subset
\delta
\}=\coprod_{\S \subset \delta}\mathcal Z_\S$$
For $\delta=\{0\}$, we put $\mathcal X_{\{0\}}=F^{ss}$. Clearly the
$X_\delta$  are disjoint  $G$-invariant subsets such that
$$ F= \coprod_{\delta} \mathcal X_\delta .$$

For any $\S,\S'$ in the same class $\delta \in \Sigma_{\ad}$, we
may define an isomorphism between the manifolds $\V(\S)$ and
$\V(\S')$ by choosing suitable representatives $V_s$ and $V_{s'}$.
This gives an isomorphism between the Hamiltonian quotients
$\mathcal Q_{\sigma_\S}$ and $\mathcal Q_{\sigma_{\S'}}$.

\

we have proved :
\begin{theorem}\label{principal}
Let $(F,G,\alpha)$ be  a general factorization problem with an
energy-complete symplectization  $\sigma$. Then we may define a
stratification
$$F= \coprod_{\delta \in \Sigma_{\ad}} \mathcal X_\delta$$
by $G$-invariant subsets defined by :
\begin{itemize}
\item $\mathcal X_{\{0\}}$ consists of the subset $F^{ss}$ of
$\sigma$-semistable elements;
\item for a non trivial class $\delta$, the stratum  $X_\delta$ is a
disjoint union
$$\mathcal X_\delta= \coprod_{\substack{\S \in (H(G)/\sim) \\
\S\subset \delta}} \mathcal Z_\S$$
where $$ Z_\S=\{f\in F \, | \, f \text{ is non $\sigma $-semistable
and } \Lambda_f = \S\}.$$
\end{itemize}
We have natural quotient maps $Z_\S \to \mathcal Q_{\sigma_\S}$,
where $\mathcal Q_{\sigma_\S}$ is the Hamiltonian quotient associated
to the factorization problem $(\V_\S,G(\S)/U(\S),\alpha_\S)$ and to
the symplectization $\alpha_\S$.

For any $\mathcal Z_\S,\mathcal Z_{\S'}$ in $\mathcal X_\delta $, the
Hamiltonian quotient  $\mathcal Q_{\sigma_\S}$ and $\mathcal
Q_{\sigma_{\S'}}$ are isomorphic.
\end{theorem}

As we will see in the last section, for the examples we have
computed, it remains that there are only  a finite number of
classes in $\Sigma_{\ad}$ which may correspond to the class of an
optimal destabilizing element, so that the number of stratum is
finite. We believe that, at least for a large class of actions,
this is  the general behavior.

\section{Linear actions}\label{section6}

We focus here  our attention on linear actions. This is  a special
case  of the previous chapter. In this case, it is possible to be
more accurate concerning the definition of the associated
factorization problem. Indeed,  it can be built as a quotient
vector subspace. Moreover the induced action is much  more
understandable.

\

So, let $\rho:G \to GL(V)$ be a linear action of a reductive group $G$
on a finite dimensional vector space $V$.

Fix a maximal compact subgroup $K$ of $G$ and an
$\ad_G$-invariant inner product of real type on $\g$. If $h$ is a
$K$-invariant Hermitian inner product on $V$, one has a standard
moment map for the $K$ action which is given by
$$\mu_0(v)=\rho^\star(-\frac{i}{2}v\otimes v^\star)$$
and any other moment map has the form
$$\mu_\tau= \mu_0-i\tau$$
with  $\tau \in iz(\mathfrak k)$. So we get a symplectization
$\sigma =(K,h,\mu_\tau)$ for the $\rho$ action. Let us remark that
in the case of a linear action, the symplectization is always
energy-complete and thus produces a well defined weight map
$\lambda^\tau  :H(G) \to \R\cup\{\infty\}$.

\

Now, for each $\xi \in i\mathfrak k$, we can decompose $V$ into
eigenspaces $V=\bigoplus_{i=1}^k V_i$ where
$\rho_\star(\xi)|_{V_i} = \xi_i id_{V_i}$ and $\xi_i$ are the
distinct eigenvalues of $\xi$. Now we have a very simple
expression for $\lambda^\tau(\xi,v)$: put
$$V_\xi^{\pm}:= \bigoplus_{\pm \xi_i>0} V_i, \ \
V_{\pm}^\xi:=\bigoplus_{\pm \xi_i\geqslant 0} V_i$$ Any $v\in V$
decomposes as $v=\sum_{i=1}^k v_i$ with $v_i \in V_i$. Then, we can
compute the map $\lambda $ in the following way~:
$$ \hspace*{-.1cm}\lambda^\tau  (v,\xi):=\lim_{t\to
+\infty}\langle\mu_\tau(\rho(e^{t\xi})v),-i\xi\rangle =
\begin{cases} +\infty
\text{ if $\exists i$ \,  s.t. $\xi_i>0$ and $v_i\not= 0$; } \\
\langle\tau,\xi\rangle \text{ otherwise}
\end{cases}$$

\

   Let $\S$ be a non trivial equivalence class of normalized
Hermitian type elements and let $\xi \in \S\cap i\k$ with
$\langle\tau,\xi\rangle<0$. Then
$$ \mathcal Z_\S=  \Bigg\{ v \in V \, | \, v\in V_-^\xi \text{ and }
\langle\tau ,\xi
\rangle=\min\limits_{\begin{array}{cc}\scriptstyle\zeta\in i\kg, \,
\|\zeta\|=1\\
\scriptstyle v\in V_-^\zeta \end{array}}\langle\tau ,\zeta \rangle
   \Bigg\}$$
The complex manifold associated to $\S$ is the complex space
   $$\V(\S)= V_-^\xi /V_\xi ^-.$$
    Let us remark that this vector space comes with a natural action of
$G(\S)=G(\xi )$ since this parabolic subgroup leaves the flag $V_\xi
^- \subset V_-^\xi$ invariant and that $U(\xi )$ acts trivially on
the quotient. So we get a well defined action $\alpha_\S$  of
$G(\S)/U(\S)$ over $ V_-^\xi /V_\xi ^-$.

We may take as a representative for the symplectization $\sigma_\S$
introduced above
(see. \ref{natsymp}) the triple
$$ (K\cap Z(\xi ), h_{|(V^-_\xi)^{\bot_h}},  \mu'= i^\star
\mu_{|(V^-_\xi)^{\bot_h}} - \langle\tau,\xi\rangle\langle i\xi ,
\cdot\rangle)$$
where $(V^-_\xi)^{\bot_h}$ denotes the orthogonal of $V^-_\xi$ in $ V_-^\xi$.

Let $v\in \mathcal Z_\S$, and let $v_0$ be the projection onto $
V_-^\xi /V_\xi ^-$.

In this framework, our general result \ref{harder}  becomes

\begin{proposition}
The vector $v_0 \in V_-^\xi /V_\xi ^-$ is $\sigma_\S$-semistable.
\end{proposition}
We give below a simple self-contained proof of this result.
\begin{proof}
    Denote by $\lambda'$ the map associated to the symplectization
$\sigma_\S$. Let $s\in i\mathfrak k\cap \mathfrak z(\xi)$, then,
$s$ admits an orthogonal decomposition as $s=\beta\xi + s^\bot$.

Assume that $\lambda '(s,v_0)<0$ so that $v_0\in {\V(\S)}^s_-$.
Using the fact that $s$ and $\xi $ commute and so are
simultaneously diagonalizable we get $v_0\in {\V(\S)}^{s^\bot}_-$
and:
\begin{align*}
\lambda '(s,v_0) &= \langle\tau ,s\rangle - \langle\xi
,s\rangle\lambda _{inf} \\
                &= \langle\tau ,s\rangle - \langle\xi ,s\rangle
\langle\tau ,\xi \rangle \\
&= \langle\tau ,s^\bot\rangle \\
&= \lambda'(s^\bot,v_0)
\end{align*}

Let $\xi_\varepsilon= \xi + \varepsilon s^\bot$ for $\varepsilon
>0$.  Using again the fact that $\xi $ and $s^\bot$ are
  simultaneously diagonalizable, it is easy to see that for
  $\varepsilon $ small enough $v\in V^{\xi _\varepsilon }_-$.
Then we get
\begin{align*}
\lambda (\frac{ \xi_\varepsilon}{\|\xi_\varepsilon\|},v)  &=
\frac{\lambda(\xi_\varepsilon,v)}{\|\xi_\varepsilon\|} \\
&= \frac{\langle\tau ,\xi \rangle + \varepsilon \langle\tau
,s^\bot\rangle}{\|\xi + \varepsilon s^\bot\|}
\end{align*}
Now we get
$$ {\frac{d}{d\varepsilon }}|_{\varepsilon =0}(\lambda (\frac{
\xi_\varepsilon}{\|\xi_\varepsilon\|},f))= \langle\tau ,s^\bot\rangle
- \langle\tau ,\xi \rangle \langle s^\bot,\xi \rangle= \lambda
'(s,v_0) <0.$$
Thus, by taking $\varepsilon $ small enough, we get a normalized
element $\frac{\xi _\varepsilon }{\|\xi _\varepsilon \|} \in
H(G)$ with $\lambda (\frac{
\xi_\varepsilon}{\|\xi_\varepsilon\|},v)<\lambda_{\rm inf}$ which
is a contradiction.
\end{proof}

We retrieve here the natural quotient maps $\V(\S)\rightarrow {\cal
Q}_{\sigma_\S}$ defined in section \ref{stratquo}. We give in the
last section examples of such linear actions, associated
stratifications and quotients maps.

\

Let us consider the example of a  linear action $\rho $ of a complex 
torus $T$ over a
  complex vector space $V$.
\

{\it We aim to show here that even in this very simple case, the 
optimal destabilizing vector $\xi $ may not be algebraic (that is may 
not lie in the subset $OPS(G)\subset \g$).}
\

We have a decomposition
$$ V= \bigoplus_{\chi \in R}V_\chi,$$
where $R\subset \Hom(T,S^1)$ and $\rho (t)|_{V\chi} = \chi(t)\id_{V\chi}$ for
all $t\in T$.

The following picture explains geometrically how we find  the optimal
destabilizing vector associated to  a symplectisation $\sigma =(K,h,\mu_\tau
)$ and to a nonsemistable vector $v$:

Using the expression of the weight map $\lambda ^\tau $, we see that 
this optimal vector as to be search in the subset
$$ \cal C=\{ s\in \g \,|\, d_e\chi(s) \leq 0, \forall \chi \in R 
\text{ s.t. } v_\chi\not= 0 \} \cap \{s\in \g \, |\, \|s\|=1\} .$$
The two planes $\cal H_{\tau _1}$ and $\cal H_{\tau _2}$  in the 
picture represents the hyperplanes defined by the equation
$$ \langle \tau ,s \rangle = \lambda ^\tau _{\rm inf} $$
for two distinct values $\tau _1, \tau _2 \in \mathfrak t$.

\

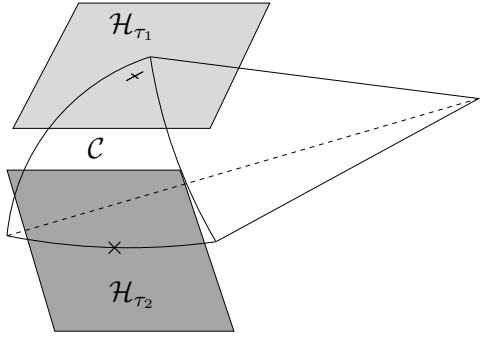
\begin{center}
\begin{picture}(0,0)%
\includegraphics{fig1.pstex}%
\end{picture}%
\setlength{\unitlength}{1973sp}%
\begingroup\makeatletter\ifx\SetFigFont\undefined%
\gdef\SetFigFont#1#2#3#4#5{%
   \reset@font\fontsize{#1}{#2pt}%
   \fontfamily{#3}\fontseries{#4}\fontshape{#5}%
   \selectfont}%
\fi\endgroup%
\begin{picture}(5949,4149)(1564,-3673)
\put(3151, 
89){\makebox(0,0)[b]{\smash{\SetFigFont{11}{13.2}{\familydefault}{\mddefault}{\updefault}{\color[rgb]{0,0,0}$\cal H_{\tau_1}$}%
}}}
\put(2701,-1486){\makebox(0,0)[b]{\smash{\SetFigFont{11}{13.2}{\familydefault}{\mddefault}{\updefault}{\color[rgb]{0,0,0}$\mathcal C$}%
}}}
\put(3151,-3286){\makebox(0,0)[b]{\smash{\SetFigFont{11}{13.2}{\familydefault}{\mddefault}{\updefault}{\color[rgb]{0,0,0}$\cal H_{\tau_2}$}%
}}}
\end{picture}
\end{center}

As $\tau $ can be choosen freely in $\z(\mathfrak t)=\mathfrak t$, 
the optimal vector may be reached at  any point of $\cal C$.

\section{Linear examples}

\subsection{Non-semistable points in the factorization problems which yield
the Grassmannians.}
{\ }\vspace{4mm}

Let $V$, $V_0$ be two Hermitian vector spaces of dimensions
$r=\dim(V)$, $r_0:=\dim(V_0)$. Consider the natural action
$\alpha_{\rm can}$ of $GL(V)$ on the space of linear morphisms
$F:=\Hom(V,V_0)$, given by $(u,f)\mapsto f\circ u^{-1}$. A moment
map for the restricted $U(V)$-action has the form
$$\mu_t(f)=\frac{i}{2}f^*\circ f-it\id_{V}\ ,\ t\in\R\ ,
$$
and the corresponding Hamiltonian quotients of $F$ are
$$Q^F_{\mu_t}=\left\{
\begin{array}{ccc}
\G r_{r}(V_0)&{\rm if}& t>0\\
\{*\}&{\rm if}& t=0\\
\emptyset&{\rm if}& t<0\ .
\end{array}\right.
$$

Fix $t>0$.  With respect to the moment map $\mu_t$ a point $f\in F$ is not
semistable if  and only if $\ker f\ne\{0\}$.  In this case, and an element
$s\in iu(V)$ destabilizes $f$ if and only if the following two conditions are
satisfied\\
\begin{itemize}
\item  $  V^-_s\subset\ker(f)$,
where $V^-_s:=\bigoplus\limits_{\begin{array}{c}\scriptstyle\lambda\in{\rm
Spec}(s)\vspace{-1mm}\\\scriptstyle \lambda<0\end{array}} V_\lambda\ ,$
\item $\lambda^s(f)=t\tr (s)<0 $.
\end{itemize}
\vspace{2mm}

This shows that the unique  normalized optimal
destabilizing   element  of $iu(V)$ is
$s_f:=-\frac{1}{\sqrt{\dim(\ker(f)}}{\rm pr}_{\ker(f)}$.

\

Let $\S\in H(gl(V))/\sim$ be the equivalence class of $s_f$. The vector
space
$V^-_\S=V^-_{s_f}$ depends only of $\S$ and the set $\mathcal Z_\S$ is
given by
$$\mathcal Z_\S=\{u \in \Hom(V,V_0) \, | \, \ker(u)=V^-_\S \}$$

   The canonically associated manifold $\V(\S)$   provided by Theorem
\ref{principal}  is
$$\V(\S):=\qmod{F}{\{u\in\Hom(V,V_0) |\
u|_{V^-_\S}=0\}}=\Hom(\qmod{V}{V^-_\S},V_0)\ .
$$
whereas the reductive quotient $G(\S)/U(\S)$ is the product
$$G_\S:=GL(V^-_\S)\times GL(V/V^-_\S)\ .
$$

The reductive group $G_\S$ acts on $\V(\S)$ in the obvious way such that
the first factor of $G_\S$ operates trivially.

The moment map $\mu'$ associated with this new action (see
\ref{natsymp}), is
$$\mu'_t:\V(\S)\ra
u(V^-_\S)\oplus u(V/V^-_\S)$$
   given by
$${\mu'_t} (\varphi)=(0,  \frac{i}{2}\varphi^*\circ \varphi
- it\id_{V/V^-_\S}) $$
and the quotient ${\cal Q}_{\sigma_\S}$ is just the Grassmannian $\G
r_{r-\dim(V^-_\S)}(V_0)$.

Therefore, applying our general  result to the factorization problem
$$(\Hom(V,V_0),GL(V),\alpha_{\rm can})$$
with the symplectization defined by $\mu_t$, $t>0$, one gets the
stratification
$$\Hom(V,V_0)=\coprod_{\rho\leq r}\Hom(V,V_0)_\rho\ $$
with
\begin{align*} \Hom(V,V_0)_\rho  := & \ \{f\in\Hom(V,V_0)|\
\rk(f)=\rho\} \\  = &\coprod_{\dim W=r-\rho } \{f\in \Hom(V,V_0)|\
\ker(f)=W\}\end{align*}
of $\Hom(V,V_0)$ and the natural quotient maps
$$\Hom(V,V_0)_\rho\map  \G r_\rho(V_0)
$$
on the strata. This is the Shatz stratification of the factorization
problem $(\Hom(V,V_0),GL(V),\alpha_{\rm can})$.

\subsection{Non-semistable points in the factorization problems which yield
the flag manifolds.}
{\ }\vspace{4mm}

Let $V_1,\dots,V_m, V=V_{m+1}$ be Hermitian vector spaces. Put
$$d_i:=\dim(V_i)\ ,\ d:=\dim(V)\ ,\  F:=\bigoplus_{i=1}^m
\Hom(V_i,V_{i+1})\ ,\ K:=\prod_{i=1}^m \U(V_i)\ ,
$$
and consider the $K$-action $\alpha_{\rm can}$ on $F$ given by
$$\alpha_{\rm can}(g_1,\dots,g_{m})(f_1, \dots,f_m)=   (g_2\circ f_1\circ
g_1^{-1},  \dots, g_m\circ f_{m-1}\circ g_{m-1}^{-1}, f_m\circ
g_{m}^{-1})\ .
$$
The general form of a moment map for the restricted  $K$-action on $F$ is
$$\mu_t(f_1,\dots,f_m)=\frac{i}{2}
\left(
\begin{array}{c}f_1^*\circ f_1 \\
f_2^*\circ f_2-f_1\circ f_1^*\\
\dots\\
f_m^*\circ f_{m}-f_{m-1}\circ f_{m-1}^*
\end{array}\right)
-i\left(\begin{array}{c} t_1\id_{V_1}\\
t_2\id_{V_2}\\\dots\\ t_m\id_{V_m}
\end{array}\right)
$$
where $t\in\R^m$.  To every $f=(f_1,\dots f_m)\in F$ we associate the
subspaces
$$W_i(f):=(f_m\circ\dots\circ f_i)(V_i)\subset V\ ,\ 1\leq i\leq m\ .$$
One obviously has $W_i\subset W_{i+1}$ and the map
$$f\mapsto (W_i(f))_{1\leq i\leq m}
$$
is constant on  orbits.
We refer to \cite{okte} for the following simple result
\begin{proposition} Suppose that   $t_i>0$, for all $1\leq i\leq m$.
\begin{enumerate}
\item Let $f\in F$.
Then the following conditions are equivalent:
\begin{enumerate}
\item $f$ is $\mu_t$-semistable
\item $f$ is $\mu_t$-stable
\item all maps $f_i$ are injective.
\end{enumerate}
\item The map
$$w:f\mapsto (W_i(f))_{1\leq i\leq m}
$$
identifies the Hamiltonian quotient $Q^F_{\mu_t}$ with the flag manifold
$$\F_{d_1,\dots d_m}(V):=\{(W_1,\dots,W_m)|\ W_1\subset
   \dots\subset W_m\subset V,\ \dim(W_i)=d_i\}\ .$$
\end{enumerate}
\end{proposition}

Fix $t=(t_1,\dots, t_m)\in\R_{>0}^m$. We  assume   $d_1\leq d_2\leq \dots \leq
d_m$, which insures that
$\F_{d_1,\dots d_m}(V)$ is non-empty. We do not require strict inequalities;
when some of the $d_i$-s coincide, the corresponding flag manifold
$\F_{d_1,\dots d_m}(V)$ can be identified with a flag manifold associated
with a smaller $m$. More precisely $\F_{d_1,\dots
d_m}(V)\simeq\F_{d_{i_1},\dots d_{i_k}}(V)$ if
$i_1<i_k<\dots <i_k$ and $\{d_1,\dots d_m\}=\{d_{i_1},\dots d_{i_k}\}$.

Suppose that $t_i>0$, for all $1\leq i\leq m$,   let
$f=(f_1,\dots, f_m)$ be a non-semistable point with respect to
$\mu_t$ and denote by $\S$ the class of its optimal destabilizing
element. The associated manifold $\V(\S)$ is
$$\V(\S):=\bigoplus_{i=1}^m \Hom(V_i^\S,V_{i+1}^\S)\ ,
$$
where $V_{m+1}^\S=V$ and $V_i^\S:=\qmod{V_i}{E_i^\S}$ with
$E_i^\S:=\ker(f_m\circ\dots \circ f_i)$ (this does not depend of the
choice of $f\in \mathcal Z_\S$).

The reductive group $G_\S$ associated with $\S$ is  the product
$$G_\S:=\prod_{i=1}^m GL(E_i^S)\times \prod_{i=1}^m
GL(V_i^\S)
$$
and the first factor operates trivially. We put  $\bar G_\S:=\prod_{i=1}^m
GL(E_i^\S)$.

The point $f_0$ of $\V(\S)$ associated with the non-semistable
point $f$ is just  $  f_0=(\bar f_1,\dots ,\bar f_m)$, where $\bar
f_i\in\Hom(V_i^\S,V_{i+1}^\S)$ is induced by $f_i$. It is easy to
see that $\bar f_i$ is injective, so the system $f_0$ defines
indeed a $(t_1,\dots,t_m)$-stable point with respect to the $\bar
G_\S$-action on $\V(\S)$. The corresponding point in the $\bar
G_\S$-quotient of $\V(\S)$ is just $(W_1(f),\dots,W_m(f)) \in
\F_{\bar d_1,\dots,\bar d_k}(V)$, where $\bar
d_i:=\rk(f_m\circ\dots \circ f_i)$.

Therefore, our general result  applied to the factorization problem
$$(\bigoplus_{i=1}^m \Hom(V_i,V_{i+1}), \prod_{i=1}^m GL(V_i),\alpha_{\rm can})
$$
with the symplectization defined by $\mu_t$   yields the
natural {\it rank}-stratification
$$\bigoplus_{i=1}^m
\Hom(V_i,V_{i+1})=\coprod_{\begin{array}{c}\scriptstyle
(\rho_1,\dots,\rho_m)\vspace{-1mm}\\
\scriptstyle\rho_1\leq\dots\leq\rho_m \vspace{-1mm}\cr\scriptstyle
0\leq\rho_i\leq d_i
\end{array}}F_{\rho_1,\dots\rho_m}\ , $$
of $\bigoplus_{i=1}^m \Hom(V_i,V_{i+1})$. The Shatz strata are
\begin{align*}
F_{\rho_1,\dots,\rho_m}:= & \{(f_1,\dots,f_m)\in \bigoplus_{i=1}^m
\Hom(V_i,V_{i+1})|\  \rk(f_m\circ\dots\circ f_i)=\rho_i\} \\
= & \hspace*{-.4cm}\coprod_{\substack{(E_1,\dots,E_m) \\ \dim
(E_i)=d_i-\rho_i}} \hspace*{-.4cm}\{(f_1,\dots,f_m)\in
\Hom(V_i,V_{i+1})| \ker(f_m\circ \dots\circ f_i)=E_i\}\ .
\end{align*}
The natural quotient maps provided by our general construction are just
the obvious maps $F_{\rho_1,\dots\rho_m}\ra \F_{\rho_1,\dots \rho_m}(V)$.

\section{Optimal destabilizing vectors in Gauge Theory}

In order to avoid the complications related to singular sheaves, we will
treat here the case when the base manifold is a complex curve $Y$. 
Another reason
for choosing this framework is the following: the natural Hamiltonian 
action of the complex gauge group on the configuration space 
associated with a {\it linear} moduli problem on a complex curve is 
{\it formally} energy complete, so it is natural to  expect  that all 
our results above can be easily generalized to this infinite 
dimensional framework.

\subsection{Holomorphic fibre bundles}

{\ }\vspace{4mm}

Let $E$ be a complex vector bundle of rang $r$ over the Hermitian curve
$(Y,g)$.  We denote by ${\cal G}$ the complex gauge group ${\cal
G}:=\Aut(E)$. Its formal Lie algebra is $A^0(\End(E))$.

The groups which play the role of  the  maximal compact subgroups    in
our gauge theoretical framework are the subgroups of the form
$${\cal
K}_h:=U(E,h)\subset {\cal G}\ ,$$
where $U(E,h)$ stands for the group of unitary automorphisms of
$E$ with respect to a Hermitian structure $h$ on $E$.

Following our general terminology developed in the finite dimensional
case, we will say that an element $s\in A^0(\End(E))$ is of Hermitian
type if there exists a Hermitian metric $h$ on $E$ such that $s\in
A^0(\Herm(E,h))$.

We are interested in the stability theory for the ${\cal
G}$-action on the space ${\cal H}(E)$ of holomorphic  structures
(semiconnections) on $E$ (see \cite{lute}).  Fixing a Hermitian
metric $h$, our moment map for the ${\cal K}_h$-action on ${\cal
H}(E)$ has the form
$$\mu({\cal E})=\Lambda_g(F_{{\cal E},h})+\frac{2\pi
i}{{Vol_g(Y)}}\frac{\deg(E)}{r} \id_E \
.
$$

One has an explicit formula for the maximal weight map $\lambda$ in this
case (see \cite{mu}).

We will need the following notation: If $a$ is an endomorphism of a vector
space
$V$, and
$\lambda\in\R$, we will put
$$V_a(\lambda):=\bigoplus_{\lambda'\leq\lambda}{\rm Eig}(a,\lambda')\ .
$$
The notation extends for endomorphisms with constant eigenvalues on
vector bundles in an obvious way.

   If ${\cal E}\in{\cal H}$ and $s\in
A^0(\Herm(E,h))$, then
$$\lambda^s({\cal E})=
\left\{\begin{array}{l}
\lambda_k\deg({\cal
E})+\sum\limits_{i=1}^{k-1}(\lambda_i-\lambda_{i+1})\deg({\cal
E}_i)-\frac{\deg({\cal E})}{r}\tr(s)\\ \hbox{if the eigenvalues}\
\lambda_1<\dots<\lambda_k
\
\hbox{of $s$ are}  \hbox{ constant and}\\ \hbox{$\EE_i:={\cal
E}_s({\lambda_i})$ are holomorphic}\\ \\
\infty \ \hbox{if not}\ .
\end{array}\right.
$$

Suppose that ${\cal E}$ is not semistable. Let
$$0=\EE_0\subset {\cal E}_1\subset {\cal E}_1\subset\dots\subset
{\cal E}_k={\cal E}
$$
be the Harder-Narasimhan filtration of ${\cal E}$ (see
\cite{hana}, \cite{bru2} for the non-algebraic case).  We recall
that this filtration is characterized by the two
conditions:\vspace{2mm}
\begin{itemize}
\item The quotients ${\cal E}_{i+1}/{\cal E}_i$ are torsion free and
semistable.\vspace{3mm} \item The slope sequence $(\mu({\cal
E}_{i+1}/{\cal E}_i))_i$ is strictly decreasing.
\end{itemize}

Put $r_i:=\rk({\cal E}_{i}/{\cal E}_{i-1})$. For any Hermitian
metric $h$ on  $E$ the optimal destabilizing element $s\in
A^0(\Herm(E,h))$ is given by the formula
$$s=\frac{1}{\sqrt{\sum\limits_{i=1}^k r_i\left[ \frac{\deg({\cal E}_{i}/{\cal
E}_{i-1})}{r_i }- \frac{\deg({\cal E})}{r}\right]^2 }}
\sum_{i=1}^k\left[\frac{\deg({\cal E})}{r}- \frac{\deg({\cal
E}_{i}/{\cal E}_{i-1})}{r_i
}\right] \id_{F_i}\ ,
$$
where $F_i$ is the $h$-orthogonal complement of ${\cal E}_{i-1}$ in
${\cal E}_{i}$.

It is not difficult to show that the holomorphic structure
$e^{st}({\cal E})$ converges to the direct sum holomorphic
structure $\bigoplus_{i=1}^k{\cal E}_{i}/{\cal E}_{i-1}$ as
$t\rightarrow \infty$. This holomorphic structure is indeed
semistable with respect to the smaller gauge group $\prod_{i=1}^k
\Aut(E_i/E_{i-1})$ and a suitable moment map.

The Shatz stratum of ${\cal E}$ is the space of all holomorphic
structures ${\cal
F}\in {\cal H}(E)$ whose Harder-Narasimhan filtration has the same
topological type
as the Harder-Narasimhan filtration  of  ${\cal E}$.

\subsection{Holomorphic pairs}

{\ }\vspace{4mm}

Let ${\cal F}_0$ be a fixed holomorphic bundle of rank $r_0$ with a fixed
Hermitian metric $h_0$  and
$E$ a complex bundle of rank $r$ on the Hermitian curve $(Y,g)$.

We are interested in the following classification problem: classify
pairs $({\cal
E},\varphi)$, where ${\cal E}$ is a holomorphic structure on $E$ and
$\varphi$ is a
holomorphic morphism $\varphi:{\cal F}_0\rightarrow {\cal E}$.  Such
a pair will be
called a holomorphic pair of type $(E,{\cal F}_0)$, and we will
denote by ${\cal
H}(E,{\cal F}_0)$ the space of such holomorphic pairs.

Our complex gauge group is ${\cal G}:=\Aut(E)$ and the role  of
the maximal compact subgroups of ${\cal G}$  are played by the
groups ${\cal K}_h:=U(E,h)$ associated with Hermitian metrics on
$E$.

   For any Hermitian metric $h$ on $E$ the moment map for the ${\cal 
K}_h$-action
on ${\cal H}(E,{\cal F}_0)$   has the form:
$$\mu({\cal E},\varphi)=\Lambda_g F_{{\cal
E},h}-\frac{i}{2}\varphi\circ\varphi^*+\frac{i}{2} t\id_E\ .
$$

{\it Suppose that   the following
topological condition holds:}
$$\mu({\cal E})\leq\tau$$
(this is
the obvious topological condition implied by the equation
$\mu({\cal E},\varphi)=0$ when one integrates its trace over $Y$).

\noindent It is well-known (\cite{bra}) that in this case a holomorphic pair
$({\cal E},\varphi)$ with $\varphi\ne 0$ is semistable with
respect to this moment map if and only if it is
$\tau:=\frac{1}{4\pi} t \Vol_g(Y)$-semistable in the following
sense:

\vspace{1.5mm}
\begin{enumerate}
\item
$\frac{\deg(\FF)}{\rk(\FF)}\leq\tau$ for all reflexive subsheaves
$\FF\subset\EE$ with $0<\rk(\FF)<r$.
\vspace{2mm}
   \item $\frac{\deg(\EE/\FF)}{\rk(\EE/\FF)}\geqslant \tau$
for all reflexive subsheaves
   $\FF\subset\EE$ with $0<\rk(\FF)<r$ and $\varphi\in H^0(\Hom(\FF_0,\FF))$.
\end{enumerate}

Note that in the case $\mu(\EE)=\tau$, the $\tau$-semistability of the pair
$(\EE,\varphi)$ is equivalent to the semistability of the bundle $\EE$. Such a
pair can be $\tau$-polystable only if $\varphi=0$.

Using the same method as in the case of bundles one obtains the
following analogue of the Harder-Narasimhan  theorem (see \cite{brte2}).
\begin{theorem}\label{vortex} Let $({\cal E},\varphi)$ be a non
$\tau$-semistable holomorphic pair of type $(E,{\cal F}_0)$ with 
$\mu({\cal E})\leq\tau$. Then there exists a unique holomorphic 
filtration
with torsion free quotients
$$0=\EE_0\subset \EE_1
\subset\dots\subset\EE_m\subset\EE_{m+1}\subset\dots\subset\EE_k=\EE
$$
of $\EE$ such that:
\begin{enumerate}
   \item The
slopes sequence satisfies:
$$\hspace{-.6cm}\frac{\deg(\EE_{1}/\EE_0)}{\rk(\EE_{1}/\EE_0)}>\dots
>\frac{\deg(\EE_{m}/\EE_{m-1})}{\rk(\EE_{m}/\EE_{m-1})}>\tau\geq
\frac{\deg(\EE_{m+1}/\EE_{m})}{\rk(\EE_{m+1}/\EE_{m })}>
\dots>\frac{\deg(\EE_{k}/\EE_{k-1})}{\rk(\EE_{k}/\EE_{k-1})}\ .$$
\item The quotients $\EE_{i+1}/\EE_i$ are semistable  for $i\ne m$.
\item One of the following properties holds:
\begin{enumerate}
\item
$$\im(\varphi)\not\subset{\cal E}_m\ ,\ \tau>
\frac{\deg(\EE_{m+1}/\EE_{m})}{\rk(\EE_{m+1}/\EE_{m })}$$
   and the pair
$(\EE_{m+1}/\EE_{m},\bar\varphi)$ is $\tau$-semistable, where $\bar\varphi$ is
the \\$\EE_{m+1}/\EE_{m}$-valued morphism induced by $\varphi$.
\item
$$\im(\varphi)\not\subset{\cal E}_m\ ,\ \tau=
\frac{\deg(\EE_{m+1}/\EE_{m})}{\rk(\EE_{m+1}/\EE_{m })}
$$
and $\EE_{m+1}/\EE_{m}$ is semistable of slope $\tau$. This implies that the
pair $(\EE_{m+1}/\EE_{m},\bar\varphi)$ is $\tau$-semistable.
\item $\im(\varphi)\subset {\cal E}_m$ and $\EE_{m+1}/\EE_{m}$ is semistable.
\end{enumerate}
\end{enumerate}

Moreover, in the cases  (b) and (c) the obtained filtration coincides with the
Harder-Narasimhan filtration of ${\cal E}$.
\end{theorem}

One can again give an explicit  formula for the maximal weight function which
corresponds to our gauge theoretical problem. The result is

$$\lambda^s({\cal E})=
\left\{\begin{array}{l} \lambda_k\deg({\cal
E})+\sum\limits_{i=1}^{k-1}(\lambda_i-\lambda_{i+1})\deg({\cal
E}_i)-\tau\tr(s)\\ \hbox{if the eigenvalues}\
\lambda_1<\dots<\lambda_k \ \hbox{of $s$ are}  \hbox{ constant,
$\EE_i:={\cal E}_s({\lambda_i)}$}\\ \hbox{are holomorphic, and
$\varphi\in
H^0(\Hom(\FF_0,\EE_s(0)))$}\ .\\
\\
\infty \ \hbox{if not}\ .
\end{array}\right.
$$

Put again $r_i:=\rk({\cal E}_{i}/{\cal E}_{i-1})$.

One can prove the following result:
\begin{theorem}
For any
Hermitian metric $h$ on $E$ the optimal destabilizing element
$s\in A^0(\Herm(E,h))$ of the holomorphic pair $(\EE,\varphi)$ is
given by
\begin{enumerate}
\item If $\im(\varphi)\subset \EE_m$ then
$$s=\frac{1}{\sqrt{\sum\limits_{i=1}^kr_i\left[
\frac{\deg({\cal E}_{i}/{\cal E}_{i-1})}{r_i }- \tau\right]^2
}} \sum\limits_{i=1}^k\left[\tau- \frac{\deg({\cal E}_{i}/{\cal
E}_{i-1})}{r_i }\right]
\id_{F_i}\ .
$$
\item If $\im(\varphi)\not\subset \EE_m$

$$s=\frac{1}{\sqrt{\sum\limits_{\begin{array}{c}\scriptstyle
i=1\vspace{-1.5mm}\\
\scriptstyle i\ne m+1\end{array}}^kr_i\left[
\frac{\deg({\cal E}_{i}/{\cal E}_{i-1})}{r_i }- \tau\right]^2
}} \sum\limits_{\begin{array}{c}\scriptstyle i=1\vspace{-1.5mm}\\
\scriptstyle i\ne m+1\end{array}}^k\left[\tau- \frac{\deg({\cal E}_{i}/{\cal
E}_{i-1})}{r_i }\right]
\id_{F_i}\ ,
$$
\end{enumerate}
where $F_i$ is the $h$-orthogonal complement of ${\cal E}_{i-1}$
in ${\cal E}_{i}$.
\end{theorem}
Note that in the second case the \mbox{$(m+1)$-th} eigenvalue of
$s$ vanishes.

\

The two pictures below explain geometrically why the optimal 
destabilizing vector is given by different formulae in the two cases 
$\im(\varphi)\subset \EE_m$,  $\im(\varphi)\not\subset \EE_m$.

The gray region represents the set $Z$ of those $\zeta\in A^0({\rm 
Herm}(E,h))$ of norm less or equal than $1$ with constant eigenvalues 
$\zeta_1<\zeta_2<\dots<\zeta_k$  such that: %
\begin{itemize}
\item the associated filtration $(\EE(\zeta_i))_i$ coincides with the 
filtration given by Theorem \ref{vortex},  and
\item $\lambda^\zeta(\EE,\varphi)<\infty$.
\end{itemize}

This region is a convex subset $Z$ of the space $\R^k$.

The second condition means that $\zeta_j\leq 0$ for those $j$ for 
which the projection of $\varphi$ on $F_j$ does not vanish.

The line in the two pictures represents the hyperplane $H\subset 
\R^k$ given by the  equation $ \zeta_k\deg({\cal 
E})+\sum\limits_{i=1}^{k-1}(\zeta_i-\zeta_{i+1})\deg({\cal
E}_i)-\tau\tr(s)=\lambda^{\min} (\EE,\varphi)$.

\

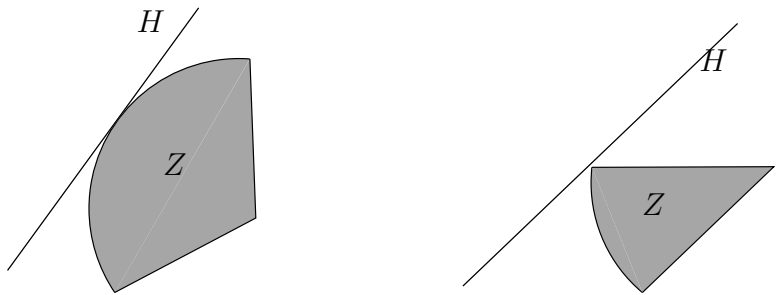
\begin{center}
\begin{tabular}{ccc}
\begin{picture}(0,0)%
\includegraphics{fig2.pstex}%
\end{picture}%
\setlength{\unitlength}{3947sp}%
\begingroup\makeatletter\ifx\SetFigFont\undefined%
\gdef\SetFigFont#1#2#3#4#5{%
   \reset@font\fontsize{#1}{#2pt}%
   \fontfamily{#3}\fontseries{#4}\fontshape{#5}%
   \selectfont}%
\fi\endgroup%
\begin{picture}(1583,1813)(1189,-1937)
\put(2101,-286){\makebox(0,0)[b]{\smash{\SetFigFont{12}{14.4}{\familydefault}{\mddefault}{\updefault}{\color[rgb]{0,0,0}$H$}%
}}}
\put(2251,-1186){\makebox(0,0)[b]{\smash{\SetFigFont{12}{14.4}{\familydefault}{\mddefault}{\updefault}{\color[rgb]{0,0,0}$Z$}%
}}}
\end{picture} &
\hspace*{2cm}&
\begin{picture}(0,0)%
\includegraphics{fig3.pstex}%
\end{picture}%
\setlength{\unitlength}{3947sp}%
\begingroup\makeatletter\ifx\SetFigFont\undefined%
\gdef\SetFigFont#1#2#3#4#5{%
   \reset@font\fontsize{#1}{#2pt}%
   \fontfamily{#3}\fontseries{#4}\fontshape{#5}%
   \selectfont}%
\fi\endgroup%
\begin{picture}(1977,1715)(589,-1314)
\put(2176, 
89){\makebox(0,0)[b]{\smash{\SetFigFont{12}{14.4}{\familydefault}{\mddefault}{\updefault}{\color[rgb]{0,0,0}$H$}%
}}}
\put(1801,-811){\makebox(0,0)[b]{\smash{\SetFigFont{12}{14.4}{\familydefault}{\mddefault}{\updefault}{\color[rgb]{0,0,0}$Z$}%
}}}
\end{picture}
\end{tabular}
\end{center}

In the first case $H$ touches $Z$  in a smooth point of its boundary 
which belongs to the interior of the spherical  region of this 
boundary, whereas in the second case $H$ touches $Z$ in a singular 
point of its boundary.

\

One can show that $e^{st}({\cal E},\varphi)$
converges either to the   object
$$(\EE_1/\EE_0,\dots,\EE_m/\EE_{m-1},
(\EE_{m+1}/\EE_{m},\bar\varphi),\EE_{m+2}/\EE_{m+1},\dots,\EE_k/\EE_{k-1})\ ,
$$
if $\im(\varphi)\not\subset\EE_m$, or to the object
$$(\EE_1/\EE_0,\dots,\EE_m/\EE_{m-1},
\EE_{m+1}/\EE_{m},\EE_{m+2}/\EE_{m+1},\dots,\EE_k/\EE_{k-1})\ ,
$$
if $\im(\varphi)\not\subset\EE_m$.

In both cases the limit object is semistable with respect to the gauge group
$\prod_{i=1}^k
\Aut(E_i/E_{i-1})$ and a suitable moment map.

Therefore our principle holds again: the optimal destabilizing
vector of a non-semistable pair gives the generalized
Harder-Narasimhan filtration and the limit object  in the direction
of this vector is a semistable object for a new moduli problem.

Theorem \ref{vortex} allows one to define a Shatz stratification on the space
of holomorphic pairs of type $(E,\FF_0)$.

Details will appear in a forthcoming article \cite{Br}.

\

\

\nocite*

\bibliographystyle{plain}
\bibliography{article}

\end{document}